\documentclass[preprint,11pt]{elsarticle}

\usepackage{siunitx}
\usepackage{booktabs}
\usepackage{amsthm}
\usepackage{subcaption}

\usepackage{amsmath,amssymb,bm}
\usepackage{accents}

\biboptions{sort}
\usepackage[dvipsnames,svgnames,table]{xcolor}

\usepackage{color}

\usepackage{tikz}
\usetikzlibrary{plotmarks}
\usetikzlibrary{positioning}
\usetikzlibrary{decorations.pathreplacing}
\usetikzlibrary{math}

\usepackage{graphicx}
\usepackage{calc,fancyvrb}
\usepackage[margin=1.in]{geometry}

\usepackage{algorithm,algpseudocode}
\usepackage{algorithmicx}
\algrenewcommand\alglinenumber[1]{\footnotesize #1:} 

\usepackage{pgfplots}
\usepackage{graphicx}
\pgfplotsset{compat=1.18}
\usepackage[bookmarks=true,colorlinks=true,linkcolor=blue]{hyperref}

\newtheorem{theorem}{Theorem}
\newtheorem{proposition}{Proposition}[section]
\newtheorem{lemma}{Lemma}[section]

\newtheorem{remark}{Remark}[section]

\numberwithin{equation}{section}

\begin{document}

\begin{frontmatter}
\title{A geometrically robust unfitted boundary algebraic equation method based on discrete potentials and local basis functions}

\author[add1,add2]{Qing~Xia}
\ead{qxia@kean.edu}

\address[add1]{Department of Mathematics, Wenzhou Kean University, Zhejiang, China, 325060.}
\address[add2]{International Frontier Interdisciplinary Research Institute, Wenzhou Kean University, Zhejiang, China, 325060.}

\begin{abstract}
We present an unfitted boundary algebraic equation (BAE) method for solving elliptic partial differential equations in complex geometries. The method employs lattice Green's functions on infinite regular grids combined with discrete potential theory to construct single and double layer potentials, which is a discrete analog to boundary integral method. Local basis functions on cut cells accommodate arbitrary boundary conditions and seamlessly integrate with the boundary algebraic equations. The difference potentials framework enables efficient treatment of nonhomogeneous terms and fast computation of layer potentials via FFT-based solvers. We establish theoretical stability and convergence through a novel interpolation operator framework. Key advantages of the developed method include: dimension reduction, geometric flexibility, mesh-independent conditioning, small-cut stability, and uniform treatment of smooth and non-smooth geometries. Numerical experiments validate accuracy and robustness across ellipses and diamonds with varying aspect ratios and sharp corners, and an application of potential flows in unbounded domains.
\end{abstract}


\end{frontmatter}



\section{Introduction}\label{sec:introduction}
The numerical solution of elliptic partial differential equations in complex geometries emains a fundamental challenge in computational mathematics. Traditional finite difference and finite element methods require mesh generation that conforms to the domain boundary, which becomes increasingly difficult and computationally expensive for intricate geometries, especially for moving geometries. This challenge has motivated the development of unfitted methods that decouple the computational mesh from the geometric boundary or material interface.

Among unfitted approaches, immersed boundary/interface methods \cite{peskin2002immersed,leveque1994immersed}, cut finite element methods \cite{hansbo2002unfitted}, immersed finite element method \cite{zhang2004immersed}, kernel free boundary integral method \cite{ying2007kernel} and fictitious domain methods \cite{glowinski1994fictitious} have gained significant attention. However, these methods often face challenges including: reduced accuracy near boundaries due to irregular stencils, modification of stencils near the boundary or interface, conditioning problems arising from small cut cells, complex implementation requirements for maintaining stability, and etc.

Boundary integral methods offer an attractive alternative by reducing the dimensionality of the problem from the domain to its boundary. These methods naturally handle complex geometries and unbounded domains while providing excellent conditioning properties. Classical boundary integral formulations require evaluation of singular integrals and sophisticated quadrature rules \cite{helsing2013solving,klockner2013quadrature}, particularly for geometries with corners or cusps.

The discrete potentials theory has been developed independently in \cite{duffin1953discrete,saltzer1958discrete}. Single and double layer discrete surface potentials in the finite difference setting are also studied in \cite{tsynkov2003definition}. The boundary algebraic equation (BAE) method, coined by Martinsson and Rodin \cite{martinsson2009boundary}, replaces continuous Green's functions with lattice Green's functions (LGFs) \cite{economou2006green,pozrikidis2014introduction} defined on infinite regular grids. The BAE method inherits the dimensional reduction of boundary integral equation method, while avoiding singular integrals entirely. The LGF remains finite even when source and target coincide. The resulting algebraic equations exhibit spectral properties analogous to Fredholm equations of the second kind, with bounded or slowly growing condition numbers under mesh refinement.

This paper presents a novel unfitted boundary algebraic equation method that bridges the gap between unfitted finite difference methods and discrete potentials theory. The key novelty lies in employing lattice Green's functions on regular infinite grids with unfitted finite difference method. We use the discrete potential theory \cite{martinsson2009boundary} to construct single and double layer potentials, with no need to evaluate singular integrals. Local polynomial basis functions on cut cells enable accurate treatment of boundary conditions and allows seamless integration with the boundary algebraic equations, while maintaining the simplicity of regular grid computations. In essence, the basis functions serve to interpolate the discrete potentials at the boundary to approximate its continuous counterpart.

The developed method inherits several advantages from both frameworks: dimensional reduction to near-boundary only problems as in boundary integral methods, geometric flexibility without mesh generation requirements, mesh-independent conditioning for appropriate formulations, stability despite arbitrarily small cut cells, uniform treatment of smooth and non-smooth geometries, no need of artificial boundary conditions for unbounded domains, and efficient solution via FFT-based solvers via the difference potentials framework \cite{ryaben2012method}.

We establish theoretical stability and convergence through an equivalent interpolation operator framework that avoids extrapolation near boundaries, following Larsson and Thom\'{e}e \cite{larsson2003partial}. The analysis reveals that the method maintains stability even in the presence of small cut cells, addressing a critical limitation of many unfitted approaches. Extensive numerical experiments validate the theoretical predictions across diverse geometric configurations including ellipses and diamonds with extreme aspect ratios, sharp corners, and unbounded domains.

The remainder of this paper develops the method systematically. Section \ref{sec:bae} introduces lattice Green's functions and constructs the single and double discrete layer potentials. Section \ref{sec:boundary_closure} develops the boundary discretization using local basis functions and derives the resulting algebraic systems. Section \ref{sec:convergence} proves stability and convergence through the interpolation operator framework. Section \ref{sec:dpm} demonstrates how difference potentials enable fast evaluation. Section \ref{sec:nonhomogeneous_equations} extends the formulation to nonhomogeneous problems in the difference potentials framework. Section \ref{sec:num} presents comprehensive numerical experiments validating the theoretical predictions and geometric robustness. Section \ref{sec:conclusion} concludes with a discussion of extensions and applications.
\section{Boundary Algebraic Equation}\label{sec:bae}

\subsection{Lattice Green's function}
Consider the difference Laplace equation defined on the infinite grids $\mathbb{Z}^2$:
\begin{align}\label{eqn:homo}
[Au](m): = 4u(m)-u(m+e_1)-u(m-e_1)-u(m+e_2)-u(m-e_2) = 0,
\end{align}
where $m\in \mathbb{Z}^2$, and $e_1=[1,0]$, $e_2=[0,1]$ are the unit vectors in $\mathbb{Z}^2$. Fourier analysis of $[Au](m)=\delta(m)$ (where $\delta(m)$ is the discrete delta function at point $m$) shows that the fundamental solution of the difference Laplace equation is given by
\begin{align}\label{eqn:lgf1}
\widetilde{G}(m) = \frac{1}{(2\pi)^2}\int_{-\pi}^{\pi}\int_{-\pi}^{\pi} \frac{\cos(\xi_1m_1+\xi_2m_2)}{4-2\cos(\xi_1)-2\cos(\xi_2)}d\xi_1\,d\xi_2,
\end{align}
which is also known as the lattice Green's function. However, the LGF in \eqref{eqn:lgf1} is strongly singular and diverges for all values of $m$, but can be normalized:
\begin{align}\label{eqn:lgf}
G(m) = \frac{1}{(2\pi)^2}\int_{-\pi}^{\pi}\int_{-\pi}^{\pi} \frac{\cos(\xi_1m_1+\xi_2m_2)-1}{4-2\cos(\xi_1)-2\cos(\xi_2)}d\xi_1\,d\xi_2,
\end{align}
The singularity $\widetilde{G}(0)$ is subtracted off to ensure convergence of the normalized integral \cite{martinsson2002asymptotic}.

The lattice Green's functions have been studied in many physics applications, such as in quantum or solid state physics \cite{economou2006green}. It can be checked that
\begin{align}
G(0,0) = 0,
\end{align}
then by the definition $[AG](0,0) = 1$ and the sign symmetry of $G$, we have
\begin{align}
G(1,0)=G(-1,0)=G(0,1)=G(0,-1)=-\frac{1}{4}.
\end{align}

In \cite{cserti2000application}, it was established (with some modifications) that
\begin{align}\label{eqn:alt-lgf}
G(m) = \frac{1}{2\pi}\int_0^\pi \frac{e^{-|m|s}\cos ny-1}{\sinh s}dy,
\end{align}
where $\cosh s = 2-\cos y$. Equation~\eqref{eqn:alt-lgf} is more numerically stable compared to \eqref{eqn:lgf}.
For $G(1,1)$, the following analytic result can be obtained
\begin{align}
G(1,1) = -\frac{1}{2\pi}\int_0^\pi \frac{(1-\cos y)^2}{\sqrt{(2-\cos y)^2-1}}dy = -\frac{1}{\pi}.
\end{align}
The following recursion formula is also established in \cite{morita1971useful} for $j\geq1$, $0<k<j$.
\begin{subequations}\label{eqn:recursion}
\begin{align}
G(j+1,k+1) &= \frac{4j}{2j+1}G(j,j)-\frac{2j-1}{2j+1}G(j-1,j-1),\\
G(j+1,k) &= 2G(j,j)-G(j,j-1),\\
G(j+1,0)&=4G(j,0)-G(j-1,0)-2G(j,1),\\
G(j+1,k)&=4G(j,k)-G(k-1,k)-G(j,k+1)-G(j,k-1).
\end{align}
\end{subequations}
Using $G(0,0),G(1,0),G(1,1)$, symmetry $G(m_1,m_2)=G(m_2,m_1)$ and the recursion formulas, in theory, we can precompute the LGF for any positive indices $(m_1,m_2)$. However, the recursion formula \eqref{eqn:recursion} is not stable numerically for large values of $|m|$.

It has been established in, for example \cite{gillman2014fast}, that as $|m|\rightarrow\infty$, LGF~\eqref{eqn:lgf} admits the following asymptotic expansion
\begin{align}
G(m) = &-\frac{1}{2\pi}\left(\log|m|+\gamma+\frac{\log 8}{2}\right)+\frac{1}{24\pi|m|^{6}}\Big(m_1^4-6m_1^2m_2^2+m_2^4\Big)\nonumber\\
& +\frac{1}{480\pi|m|^{12}}\Big(43m_1^8-772m_1^6m_2^2+1570m_1^4m_2^4\nonumber\\
&\quad \quad \quad \quad \quad \quad-722m_1^2m_2^6+43m_2^8\Big)+\mathcal{O}\left(\frac{1}{|m|^6}\right),
\end{align}
where $\gamma$ is the Euler constant.

In our computation, for small $|m|$, we evaluate the integral \eqref{eqn:alt-lgf} directly; for large $|m|$, we use the following equivalent expansion in \cite{martinsson2002asymptotic} for better accuracy which avoids large number multiplication or division:
\begin{align}
G(m) = &-\frac{1}{2\pi}\left(\log|m|+\gamma+\frac{\log 8}{2}\right)+\frac{\cos(4\theta)}{24\pi|m|^{2}}\nonumber\\
& +\frac{25\cos(8\theta)+18\cos(4\theta)}{480\pi|m|^{4}}+\frac{490\cos(12\theta)+459\cos(8\theta)}{2016\pi|m|^6}+\mathcal{O}\left(\frac{1}{|m|^8}\right),
\end{align}
where $\theta=atan2(y,x)$.

\begin{remark}
The values of $G(m_1,m_2)$ can be evaluated with machine precision accuracy for sufficiently large $|m|$. One distinction from the Green's function for the continuous Laplace equations is that no singularity appears even when the source and target points coincide.
\end{remark}

\begin{remark}
In this work, we focus on the lattice Green's functions on simple square lattices. Other lattices such as triangular or honeycomb lattices \cite{pozrikidis2014introduction} or high order versions \cite{gabbard2024lattice} can be treated similarly and can be integrated into the unfitted method developed in this work seamlessly.
\end{remark}

\subsection{Single and double-layer potentials for homogeneous difference equations}
Assume now we have a bounded domain $\Omega\subset\mathbb{R}^2$.
We first embed the domain $\Omega$ into an infinite mesh $(\mathbb{Z}h)^2$ with grid spacing $h$, and denote $M^+$ as all points inside domain $\Omega\cup\Gamma$, and $M^-:=(\mathbb{Z}h)^2\backslash M^+$ for all points outside $\Omega$. In this work, we focus on second-order central finite difference discretization with the five-point stencil. Point sets $N^{\pm}$ are defined as follows
\begin{align}
N^\pm = \left\{(x_j,y_k),(x_{j\pm1},y_k),(x_j,y_{k\pm1})\,\vert \, (x_j,y_k)\in M^\pm\right\}.
\end{align}
The point sets $N^\pm$ intersect at $\gamma:=N^+\cap N^-$, which we will divide into 2 subsets based on their location inside or outside of the domain $\Omega$: $\gamma_+:=\gamma\cap M^+$ and $\gamma_-:=\gamma\backslash \gamma_+$.

\begin{figure}[htbp]
\centering
\includegraphics[width=0.4\textwidth]{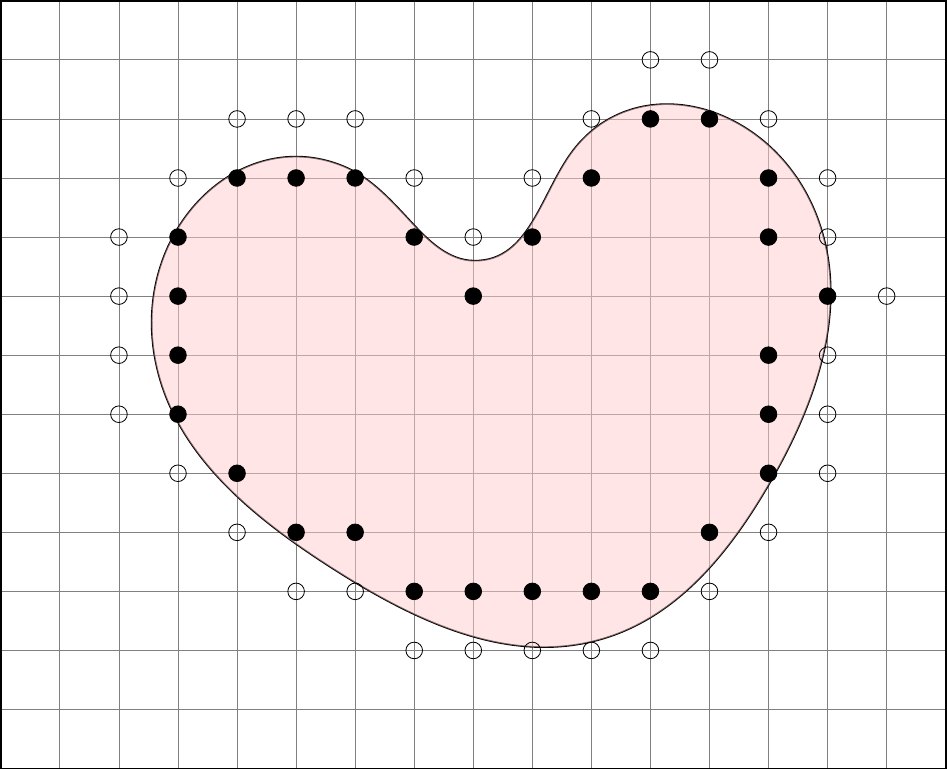}
\caption{Example of set $\gamma$ ($\gamma_+$: solid point; $\gamma_-$: circles)}\label{fig:gamma}
\end{figure}

With the point sets defined, we are ready to introduce the layer potentials.
The single-layer kernel is
\begin{align}
S(m,n) = G(m-n), \quad m,n\in \mathbb{Z}^2,
\end{align}
and the single-layer potential is
\begin{align}\label{eqn:single}
u(m) = \sum_{x_n\in{\gamma}^-}S(m,n)q(n),\quad x_m\in N^+,
\end{align}
where the discrete densities $q(n)$ are defined on the exterior point set $\gamma_-$.

\begin{proposition}
The single layer formulation \eqref{eqn:single} satisfies the difference equation \eqref{eqn:homo} inside $M^+$. 
\end{proposition}

\begin{proof}
For $x_m\in M^+$,
\begin{align}
A u_m =& A\sum_{x_n\in\gamma_-}S(m,n)q(n)=\sum_{x_n\in\gamma_-}\Big[AS(m,n)\Big]q(n)\nonumber\\
=&\sum_{x_n\in\gamma_-}\delta(m,n)q(n)=0,
\end{align}
as $m$ and $n$ belong to disjoint sets.
\end{proof}

\begin{remark}
The grid spacing $h$ is irrelevant and does not come into play as the lattice Green's function is the fundamental solution for the Laplace equation with homogeneous right hand side.
\end{remark}

Similarly, the double-layer kernel is defined as
\begin{align}\label{eqn:double}
D(m,n) = \sum_{k\in \mathbb{D}_n}G(m-n)-G(m-k), \quad m,n\in\mathbb{Z}^2,
\end{align}
where $k\in \mathbb{D}_n$ is an exterior node connected to the source point $n$ (see \cite{martinsson2009boundary} for more details).

\begin{remark}
Here we assume that the geometry is nice enough that each point in $\gamma_-$ has at least one exterior connected point in $M^-\backslash \gamma_-$.
\end{remark}

The double-layer potential is defined as
\begin{align}
u(m) = \sum_{x_n\in{\gamma}^-}D(m,n)q(n),\quad x_m\in N^+,
\end{align}
where the discrete densities are also defined on the exterior boundary point $\gamma_-$, which differs from the definition in \cite{tsynkov2003definition}.

\begin{proposition}
The double-layer formulation \eqref{eqn:double} satisfies the difference equation \eqref{eqn:homo} inside $M^+$. 
\end{proposition}
The proof is similar to the single-layer case.

Restricting our attention to the source points ${\gamma}_-$ and the target points ${\gamma}_+$, we introduce the following operators $S_\pm$ and $D_\pm$ such that
\begin{subequations}
\begin{align}
u_{{\gamma}_-} &= S_-q_s,\quad u_{{\gamma}_+} = S_+q_s,\\
u_{{\gamma}_-} &= D_-q_d,\quad u_{{\gamma}_+} = D_+q_d,
\end{align}
\end{subequations}
where $q_s,q_d$ are unknown single and double-layer densities defined on point set $\gamma_-$. Here, $S_-,D_-$ are of dimension $|{\gamma}_-|\times |{\gamma}_-|$ and $S_+,D_+$ are of dimension $|{\gamma}_+|\times |{\gamma}_-|$. 

If $\gamma_-$ happens to be the discrete grid boundary of the domain $\Omega$, as is the case for the fitted boundary algebraic equation in \cite{martinsson2009boundary}, one can solve one of the boundary equations
\begin{subequations}
\begin{align}
u_{{\gamma}_-} &= S_-q_s,\\
u_{{\gamma}_-} &= D_-q_d,
\end{align}
\end{subequations}
for the unknown density $q_s$ or $q_d$.

Once the density $q_s$ or $q_d$ is obtained, the solution to the homogeneous difference equation can be retrieved at any interior grid points inside the domain $\Omega$.

Often the case is $\gamma_-$ lies outside the domain $\Omega$ completely, then by eliminating the unknown density $q_s$ or $q_d$, we obtain a linear relation between values in $\gamma_+$ and $\gamma_-$:
\begin{subequations}
\begin{align}
u_{{\gamma}_+} &= S_+S^{-1}_{-}u_{{\gamma}_-},\\
u_{{\gamma}_+} &= D_+D^{-1}_{-}u_{{\gamma}_-},
\end{align}
\end{subequations}
both giving an effective $|\gamma_+|$ number of equations, hinting the need of additional $|\gamma_-|$ number of equations from the boundary conditions, so as to balance the number of unknowns and the number of equations.

\begin{remark}
The difference potentials method constructs the relation between $u_{\gamma_+}$ and $u_{\gamma_-}$ using the auxiliary domain and requires no explicit knowledge of the related Green's function. The lattice Green's function approach is only valid for linear and constant coefficient PDEs, and provides a computationally less expensive venue for constructing such linear relations.
\end{remark}

\section{Boundary closure}\label{sec:boundary_closure}
The lattice Green's function obtained in Sec.~\ref{sec:bae} corresponds to the second-order central finite difference schemes for Laplace equation $\Delta u = 0$. We will discuss how to handle different types of boundary conditions with second-order accuracy in this section.

\subsection{Dirichlet BC}

We first look at the Dirichlet boundary condition (BC):
\begin{align}
u(x) = g(x),\quad x\in\Gamma,
\end{align}
where $\Gamma=\partial\Omega$ is of arbitrary shape and is the boundary of the domain $\Omega$.

Following \cite{xia2023local,banks2016galerkin}, we introduce the following local basis functions defined at vertex $(x_j,y_k)$ in the infinite mesh $(\mathbb{Z}h)^2$ with grid spacing $h$:
\begin{align}
\phi_{jk}(x,y) = \phi\left(\frac{x-x_j}{h}\right)\phi\left(\frac{y-y_k}{h}\right),
\end{align}
where 
\begin{align}
\phi(\xi) = \left\{
\begin{array}{lr}
1+\xi, & -1\leq \xi \leq 0, \\
1-\xi, & 0\leq \xi \leq 1, \\
0, & \mbox{elsewhere}.
\end{array}
\right.
\end{align}
The basis functions define standard P1 element in 2D.

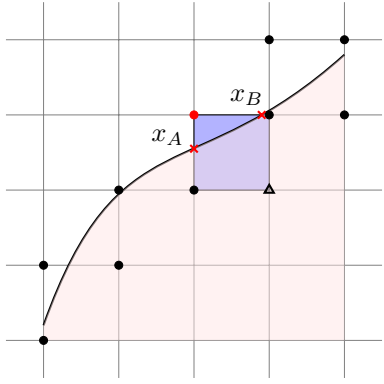
\begin{figure}[htbp]
\centering
\begin{tikzpicture}[scale=1]
\draw[step=1cm,gray,very thin] (-0.5,-0.5) grid (4.5,4.5);
\filldraw[fill=blue!30] (2,2) rectangle (3,3);
\draw[thick] (0,0.2) .. controls (1,3) and (2,2) .. (4,3.8);
\fill[red!10, opacity=0.5] 
    (0,0.2) .. controls (1,3) and (2,2) .. (4,3.8) -- 
    (4,0) -- 
    (0,0) -- 
    cycle;
\filldraw (0,1) circle (1.5pt);
\filldraw (1,1) circle (1.5pt);
\filldraw (1,2) circle (1.5pt);
\filldraw (2,2) circle (1.5pt);
\filldraw[red] (2,3) circle (1.5pt);
\draw[mark size=+2pt,thick,red] plot[mark=x] coordinates {(2,2.55)};
\filldraw (3,3) circle (1.5pt);
\filldraw (3,4) circle (1.5pt);
\filldraw (0,0) circle (1.5pt);
\filldraw (4,3) circle (1.5pt);
\filldraw (4,4) circle (1.5pt);
\draw[mark size=+2pt,thick,red] plot[mark=x] coordinates {(2.9,3)};
\draw[mark size=+2pt,thick] plot[mark=triangle] coordinates {(3,2)};
\node[left] at (2,2.7) {\small $x_{A}$};
\node[above] at (2.7,3) {\small $x_{B}$};
\end{tikzpicture}
\caption{An example of boundary points}
\label{fig:example_cut_cell}
\end{figure}

From the definition, the basis function is local in the sense that it is zero outside a two by two grid box centered at $(x_j,y_k)$ (including the box edge). This implies that using the basis function defined at points in $\gamma$ set is sufficient to discretize the Dirichlet BC, as the contribution from such as the triangle point in Fig~\ref{fig:example_cut_cell} will be zero for the two intersection points (denoted by red crosses) of boundary and the grid lines in Fig~\ref{fig:example_cut_cell}.

Also, from Fig~\ref{fig:example_cut_cell}, we can see that there are two intersection points corresponding to the red point in $\gamma_-$. As discussed in Sec.~\ref{sec:bae}, we need an additional $|\gamma_-|$ number of equations from the boundary conditions for balance of the numbers of unknowns and equations, which implies that we can choose one intersection point for each $\gamma_-$ point. In practice, we choose the one ($x_A$) closer to the point in $\gamma_-$ as a rule of thumb, although the further one $(x_B)$ will work as well.

We will denote the selected $|\gamma_-|$ number of intersection points as $x_i$, and the Dirichlet boundary condition is discretized as
\begin{align}\label{bc:dirichlet}
\sum_{x_{jk}\in\gamma} u_{jk}\phi_{jk}(x_i) = g(x_i),
\end{align}
which can be expressed in a matrix-vector form
\begin{align}
\Phi u_{\gamma} = g(x_i),
\end{align}
The matrix $\Phi$ is a sparse and is of size $|\gamma_-|\times|\gamma|$ due to the locality of the basis functions and can be further decomposed into two sub-matrices
\begin{align}
(\Phi_+ \ \Phi_-)[u_{\gamma_+}\ u_{\gamma_-}]^T=g(x_i),
\end{align}
where $\Phi_+$ is sparse and almost diagonal with size $|\gamma_-|\times|\gamma_+|$ and $\Phi_-$ is square, sparse, diagonal, and of size $|\gamma_-|\times|\gamma_-|$. The diagonal entry value is $\Phi_-$ is evaluated as the basis function at the intersection point centered at point in $\gamma_-$, which is larger if the intersection point is closer to the point. This also justifies the choice of the closer intersection point.

\begin{remark}
The local truncation error of the interpolation at the intersection point is in the order of $\mathcal{O}(h^2)$ due to the bilinear basis functions.
\end{remark}

We will use the double layer formulation as an example and comment on the difference between the single and double layer formulation.

In the double layer formulation, the boundary algebraic equation together with the discretization of boundary conditions give the following square system
\begin{subequations}
\begin{align}
u_{{\gamma}_+} = D_+D^{-1}_{-}u_{{\gamma}_-},\\
\Phi_+u_{\gamma_+} +\Phi_-u_{\gamma_-}=g(x_i),
\end{align}
\end{subequations}
with unknowns $u_{\gamma_+}$ and $u_{\gamma_-}$.
The above linear system can be reduced to a slightly small-sized linear system with unknown $u_{\gamma_-}$ only:
\begin{align}\label{eqn:double_bae}
A_du_{\gamma_-}:=\left(\Phi_+D_+D^{-1}_{-} +\Phi_-\right)u_{\gamma_-}=g(x_i).
\end{align}
The matrix $A_d=\Phi_+D_+D^{-1}_{-} +\Phi_-$ can also be regarded as a Schur complement to the original linear system, or a right preconditioned matrix when the inverse of $D_-$ is obtainable. However, it might never be wise to invert a dense matrix, especially in higher dimensions. In 2D, we compute $A_d$ numerically anyway to illustrate its spectral and conditioning property so that we can utilize this info in higher dimensions.
Once we solve the boundary algebraic equation from double layer potential \eqref{eqn:double_bae}, we will obtain $u_{\gamma_+}$ and $u_{\gamma_-}$, or equivalently $u_{\gamma}$. From $u_{\gamma_-}=D_-q_d$, we will be able to solve for the density $q_d$. 

Another approach is to solve for the double-layer density $q_d$ directly. Based on the following formulation, 
\begin{subequations}
\begin{align}
u_{\gamma_+} &= D_+ q_d,\\
u_{\gamma_-} &= D_- q_d,\\
\Phi_+u_{\gamma_+} +\Phi_-u_{\gamma_-}&=g(x_i),
\end{align}
\end{subequations}
we can reformulate the above system into a concise one:
\begin{align}\label{eqn:double_reformulate}
M_dq_d:=(\Phi_+D_+ +\Phi_-D_-)q_d = g(x_i),
\end{align}
where no inversion of matrices is needed to solve for $q_d$. Since $M_d$ is merely a convex linear combination of $D_+$ and $D_-$, we should expect similar spectral properties between $D_-$ and $M_d$, which is also confirmed in the numerical section.

With $q_d$ obtained, any interior value inside the domain can be retrieved by the double-layer formulation:
\begin{align}
u(m) &= \sum_{x_n\in\gamma_-}D(m,n)q_d(n),\label{eqn:double_convolution}
\end{align}
This approach can be applied to bounded and unbounded domains likewise and any solution at any point inside the domain can be obtained. In particular, no artificial boundary condition is needed in the unbounded case. Besides, fast summation techniques based on FMM or FFT are favorable, especially in higher dimensions.

\begin{remark}
When the single-layer kernel is used, the Schur complement becomes:
\begin{align}\label{eqn:single_bae}
A_su_{\gamma_-}:=\left(\Phi_+S_+S^{-1}_{-} +\Phi_-\right)u_{\gamma_-}=g(x_i).
\end{align}
The single-layer formulation would similarly give the no-inversion matrix:
\begin{align}\label{eqn:single_reformulate}
M_sq_s:=(\Phi_+S_+ +\Phi_-S_-)q_s=g(x_i),
\end{align} 
and $M_s$ should inherit similar spectral properties as $S_-$, which is also confirmed in the numerical results.

Hence, we can solve for the density $q_s$ directly with \eqref{eqn:single_reformulate} if we don't involve matrix inversion. The approximated solution is viable through the following convolution:  
\begin{align}
u(m) &= \sum_{x_n\in\gamma_-}S(m,n)q_s(n).\label{eqn:single_convolution}
\end{align}
\end{remark}

\begin{remark}
$M_s$ and $M_d$ are linear combinations of exterior potentials and interior potentials, evaluated at the continuous boundary, which aims to approximate their continuous counterpart but avoids evaluation of singular integrals.
\end{remark}

\begin{remark}

It is established in \cite{martinsson2009boundary} that, under proper assumptions,
\begin{align}
||D_-||_2\leq C\log(N),\quad ||S_-||_2\leq CN\log(N)
\end{align}
where $C$ is a constant independent of mesh size $h$. 

However, $D_-$ might be rank deficient in non-convex domains. It is possible that a point in $\gamma_-$ admits no exterior connection point and the definition in \eqref{eqn:double} fails. So we will require that the domain should have nice properties such that $D_-$ is non-singular if one seeks to use the double layer formulation. On the other hand, the definition of $S_-$ is always valid, ensuring the robustness of the single layer formulation, regardless of the shape or concavity of the geometry, as illustrated in the unbounded domain test case.

\end{remark}

\subsection{Robin BC}

Now assume the boundary condition is 
\begin{align}
\alpha \frac{\partial u(x)}{\partial n}+\beta u(x) = g(x),\quad x\in\Gamma,
\end{align}
where $n$ is the unit outward normal vector.
With the Robin boundary condition, the strong form \eqref{bc:dirichlet} discretized with bilinear basis functions is not sufficient for second order accuracy when normal derivatives appear in the boundary conditions. As studied in \cite{xia2023local}, the order of accuracy for Robin BC is $\mathcal{O}(h^{p})$ when degree $p$ polynomials are used as basis functions, hence we introduce local degree 2 Lagrange polynomials for Robin boundary conditions where second order accuracy can be achieved.

\begin{figure}[htbp]
    \centering
    \begin{subfigure}{0.45\textwidth}
        \centering
        \begin{tikzpicture}[scale=1]
		\draw[step=1cm,gray,very thin] (-0.5,-0.5) grid (4.5,4.5);
		\filldraw[fill=blue!30] (2,1) rectangle (3,2);
		\filldraw[fill=blue!30] (3,1) rectangle (4,2);
		\filldraw[fill=blue!30] (2,2) rectangle (3,3);
		\filldraw[fill=blue!30] (3,2) rectangle (4,3);
		\draw[thick] (0,0.2) .. controls (1,3) and (2,2) .. (4,3.8);
		\fill[red!10, opacity=0.5] 
	    (0,0.2) .. controls (1,3) and (2,2) .. (4,3.8) -- 
	    (4,0) -- 
	    (0,0) -- 
	    cycle;
		\filldraw (0,1) circle (1.5pt);
		\filldraw (1,1) circle (1.5pt);
		\filldraw (1,2) circle (1.5pt);
		\filldraw (2,2) circle (1.5pt);
		\filldraw[red] (2,3) circle (1.5pt);
		\draw (2,1) circle (1.5pt);
		\draw (3,1) circle (1.5pt);
		\draw (4,1) circle (1.5pt);
		\draw (3,2) circle (1.5pt);
		\draw (4,2) circle (1.5pt);
		\draw[mark size=+2pt,thick,red] plot[mark=x] coordinates {(2,2.55)};
		\filldraw (3,3) circle (1.5pt);
		\filldraw (3,4) circle (1.5pt);
		\filldraw (0,0) circle (1.5pt);
		\filldraw (4,3) circle (1.5pt);
		\filldraw (4,4) circle (1.5pt);
		\end{tikzpicture}
        \caption{}
        \label{fig:O3_basis}
    \end{subfigure}
    \quad
    \begin{subfigure}{0.45\textwidth}
        \centering
        \begin{tikzpicture}[scale=1]
		\draw[step=1cm,gray,very thin] (-0.5,-0.5) grid (4.5,4.5);
		\filldraw[fill=blue!30] (2,1) rectangle (3,2);
		\filldraw[fill=blue!30] (3,1) rectangle (4,2);
		\filldraw[fill=blue!30] (2,2) rectangle (3,3);
		\filldraw[fill=blue!30] (3,2) rectangle (4,3);
		\draw[thick] (0,0.5) .. controls (4.4,1) and (4.4,3) .. (0,3.5);
		\fill[red!10,opacity=0.5] (0,0.5) .. controls (4.4,1) and (4.4,3) .. (0,3.5) -- (0,0.5) -- cycle;
		\filldraw (0,1) circle (1.5pt);
		\filldraw (1,1) circle (1.5pt);
		\filldraw (1,0) circle (1.5pt);
		\draw (2,2) circle (1.5pt);
		\filldraw (2,3) circle (1.5pt);
		\filldraw (2,1) circle (1.5pt);
		\filldraw (3,1) circle (1.5pt);
		\filldraw (3,2) circle (1.5pt);
		\filldraw (2,0) circle (1.5pt);
		\filldraw (2,4) circle (1.5pt);
		\filldraw (1,3) circle (1.5pt);
		\filldraw (1,4) circle (1.5pt);
		\filldraw (0,3) circle (1.5pt);
		\filldraw (0,4) circle (1.5pt);
		\filldraw[red] (4,2) circle (1.5pt);
		\draw[mark size=+2pt,thick,red] plot[mark=x] coordinates {(3.3,2)};
		\filldraw (3,3) circle (1.5pt);
		\filldraw (0,0) circle (1.5pt);
		\draw (4,3) circle (1.5pt);
		\draw (4,1) circle (1.5pt);
		\end{tikzpicture}
        \caption{}
        \label{fig:additional}
    \end{subfigure}
    \caption{Support cells and extra points for the intersection points}\label{fig:support_cells}
\end{figure}
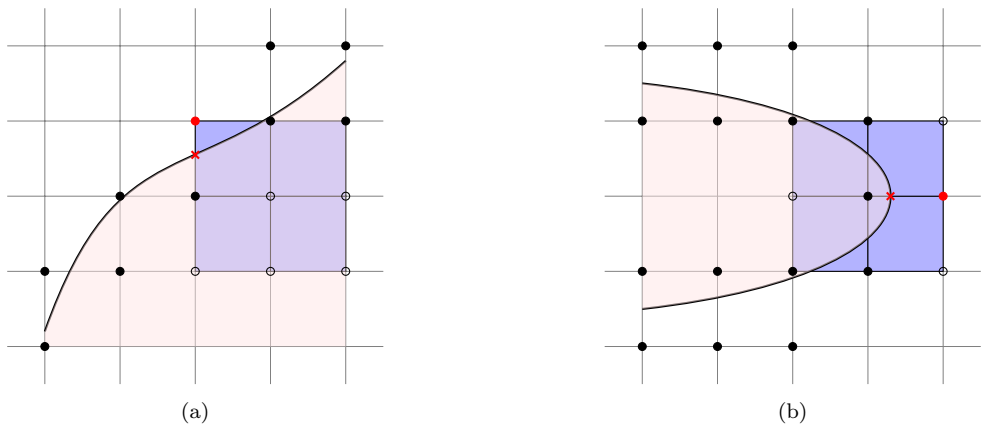

For each point in $\gamma_-$, we determine a corresponding intersection point $x_i$ between mesh line and the boundary first and look for a $3\times3$ local grid (see Fig~\ref{fig:O3_basis} and \ref{fig:additional}) that supports this intersection point. As the set $\gamma$ is not sufficient to form such $3\times3$ grid, additional points either from $M^+$ or $M^-$ will be needed. The general rule of thumb is to use as many points in $M^+$ as possible. Thus it is always possible construct such $3\times3$ grids. Including the additional points needed to form the $3\times3$ grids, we introduce the augmented sets $\tilde{\gamma}_\pm$ that incorporates the additional points depending on their location in $M_+$ or $M_-$. As also discussed in Sec.~\ref{sec:bae}, we can relate $u_{\tilde{\gamma}_+}$ and $u_{\tilde{\gamma}_-}$ using the single/double layer potentials. The number of points is irrelevant as we can put as many exterior source points as one would like. The interior target points will not alter the structure of the linear system. This will only change the number of equations for the source-target interactions.

Once we found such $3\times3$ grids, we will be able to define the local Lagrange basis functions
\begin{align}
\phi^{(2)}_{jk} = \phi^{(2)}_j(x)\phi^{(2)}_k(y),\quad j,k=1,2,3,
\end{align}
where the barycentric form is
\begin{align}
\phi^{(2)}_j(x) = \ell(x)\frac{\omega_j}{x-x_j},\quad j=1,2,3,
\end{align}
where
\begin{align}
\ell(x) = \prod_{i=1}^3(x-x_i),\quad \omega_j = \prod_{i\neq j} \frac{1}{x_j-x_i}=\frac{1}{h^2},
\end{align}
where $h$ is the uniform grid size. In $y$ direction, $\phi^{(2)}_k(y)$ is constructed similarly.

With basis function defined, we can discretize the boundary condition at the intersection point $x_i$ as
\begin{align}
\sum_{x_{jk}\in \tilde{\gamma}} u_{jk}\left(\alpha \nabla \phi^{(2)}_{jk}(x_i) \cdot n(x_i)+\beta \phi^{(2)}_{jk}\right) = g(x_i),
\end{align}
which can be also expressed into a matrix vector form $\Phi u_{\tilde{\gamma}}=g$ where $\Phi$ is still a sparse matrix with a slightly larger band of nonzero entries. It can also be decomposed into sub-matrices and sub-vectors:
\begin{align}
\left(\Phi_+ \ \Phi_{-} \ \Phi'_{-}\right)[u_{\tilde{\gamma}_+} u_{{\gamma}_-} u_\eta]^T = g(x_i)
\end{align}
where $\tilde{\gamma}_+$ denotes the $\gamma_-$ and additional interior points, and $\eta$ denotes the additional exterior points.

The differences from the Dirichlet BC include the number of unknowns and a division by $h$ in the gradient part in the coefficient matrix. The collocation points are exactly the same, which leaves the additional exterior point set $\eta$ unaccounted for. For these additional exterior points in $\eta$, we use extrapolations, which is denoted as
\begin{align}
\left(R_+ \ R_- \ I\right)[u_{\tilde{\gamma}_+} u_{{\gamma}_-} u_\eta]^T = g(x_i).
\end{align}
The exact form of $R$ depends on the geometry and might not appear for all geometries.

Consequently, the number of unknowns and the number of equations are balanced, and the coupled system becomes
\begin{subequations}
\begin{align}
u_{\tilde{\gamma}_+} &= D_+D^{-1}_{-}u_{{\gamma}_-},\\
\Phi_+u_{\tilde{\gamma}_+} +\Phi_-u_{{\gamma_-}}+\Phi'_-u_{\eta}&=g(x_i),\\
R_+u_{\tilde{\gamma}_+}+R_-u_{{\gamma}_-}+u_\eta &= 0,
\end{align}
\end{subequations}
which simplifies to the following system:
\begin{align}
A_du_{{\gamma}_-}:=\left[\Phi_+D_+D^{-1}_{-} +\Phi_--\Phi'_-(R_+D_+D^{-1}_{-}+R_-)\right]u_{{\gamma}_-}=g(x_i).
\end{align}
This is a square system where we can solve for the unknown density $u_{\gamma_-}$. If we remove all inverse matrices, we can solve for $q_d$ directly:
\begin{align}
M_dq_d:=\left[\Phi_+D_+ +\Phi_-D_{-}-\Phi'_-(R_+D_++R_-D_{-})\right]q_d=g(x_i).
\end{align}
Once the double layer density $q_d$ is obtained, we can employ the convolution \eqref{eqn:double_convolution} to obtain the solution everywhere inside the domain whether the domain $\Omega$ is bounded or unbounded.

\begin{remark}
When single layer kernel is used, the linear system we need to solve simply becomes
\begin{align}
A_su_{{\gamma}_-}:=\left[\Phi_+S_+S^{-1}_{-} +\Phi_--\Phi'_-(R_+S_+S^{-1}_{-}+R_-)\right]u_{{\gamma}_-}=g(x_i).
\end{align}
or when no inverse of matrices is used, we have
\begin{align}
M_sq_s:=\left[\Phi_+S_+ +\Phi_-S_{-}-\Phi'_-(R_+S_++R_-S_{-})\right]q_s=g(x_i).
\end{align}
Once the single layer density $q_s$ is obtained, we can use the convolution~\eqref{eqn:single_convolution} as well.
\end{remark}
\section{Convergence}\label{sec:convergence}

\begin{figure}[htbp]
\centering
\includegraphics[width=0.3\textwidth]{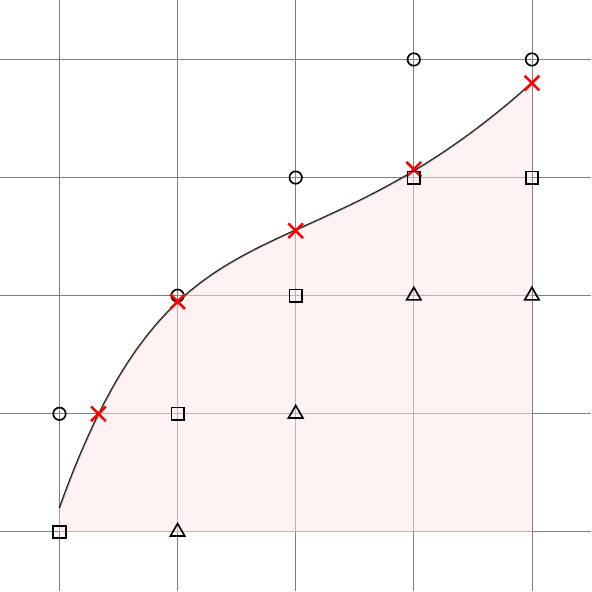}
\caption{Four sets: (1) circle $\gamma_-$; (2) red cross $\Gamma_h$; (3) square $\gamma_+$ and (4) triangle $\omega$; Shaded area denotes the interior.}\label{fig:four_sets}
\end{figure}

The solution we constructed via the boundary algebraic equation satisfies the following linear systems defined on $N^+$ and Dirichlet BC at the boundary points $x_\Gamma$ (red cross in Figure~\ref{fig:four_sets})
\begin{subequations}\label{eqn:equi_sys}
\begin{align}
\Delta_h u_m = 0, \quad x_m\in M^+,\label{eqn:laplace}\\
\alpha_m u_{m'} + (1-\alpha_m)u_{m}=g(x_\Gamma)+\tau_h,\quad x_m\in \gamma_+,\label{eqn:linear_interp}
\end{align}
\end{subequations}
where $x_{m'}$ is a point in $\gamma_-$ connected to point $x_m\in\gamma_+$ (not necessarily unique), and $x_\Gamma$ is the intersection point between the boundary $\Gamma$ and the mesh line between $x_{m'}$ and $x_m$.  The term $\tau_h$ is the truncation error in the order of $O(h^2)$. The constant $\alpha_m$ satisfies $0\leq\alpha_m< 1$. The constant $\alpha_m$ can be zero since a node on the boundary counts toward $\gamma_+$ but not $\gamma_-$.

We will categorize 4 types of point sets close to the boundary: (1) $\gamma_-$ denotes the exterior points in $\gamma$; (2) $\Gamma_h$ the boundary points; (3) $\gamma_+$ the interior points in $\gamma$; and (4) $\omega$ the extra points need to formulate the 5-point stencil at points in $\gamma_+$. Point sets $\gamma_\pm$ are used for interpolations at intersection points $\Gamma_h$. 

Next, we find an equivalent linear system from \eqref{eqn:equi_sys}
\begin{align}
Au_{j,k}=f_{j,k}+\tau_h
\end{align} 
for $x_{jk}\in M^+$, since one can replace $u_{m'}$ appeared in \eqref{eqn:laplace} with a variant form of boundary discretization from \eqref{eqn:linear_interp}
\begin{align}\label{eqn:alt_eqn}
u_{m'} = \frac{1}{\alpha_m}u(x)+\left(1-\frac{1}{\alpha_A}\right)u_{m}+\frac{1}{\alpha_m}\tau_h
\end{align}

Inside set $M^+\backslash\gamma_+$, it is easy to see that the local truncation error is 0, since the construction satisfies the difference equation exactly. For points in $\gamma_+$, we have the following conclusion.

\begin{lemma}
The local truncation error at points in $\gamma_+$ is second order accurate.
\end{lemma}

\begin{proof}
One can use the fact that the linear interpolation \eqref{eqn:linear_interp} gives second order accuracy $\tau_h\sim O(h^2)$ and at every node $x_{jk}\in\gamma_+$, we can leverage the 5-point stencil
\begin{align}
4u_{j,k}=u_{j+1,k}+u_{j-1,k}+u_{j,k+1}+u_{j,k-1}
\end{align}
where at most 3 points on the right hand side can lie in $\gamma_-$ and at least 1 point will lie in $\gamma_+$. For those lie in $\gamma_-$, we will use \eqref{eqn:alt_eqn} to eliminate $u_{\gamma_-}$, which means we will have $O(h^2)$ for the local truncation errors.
\end{proof}

Next, we prove the stability.
\begin{lemma} \label{lem:stability}
The following stability estimate holds for \eqref{eqn:equi_sys}
\begin{align}\label{eqn:ineq}
|u|_{\gamma_+} \leq \beta_1|u|_\omega+\beta_2|\tau|_{\Gamma_h}+\beta_3|\Delta_h u|_{\gamma_+},
\end{align}
where $\beta_1$ and $\beta_2$ are constants that depend on the geometry, and $\beta_1<1$.
\end{lemma}

\begin{proof}
For any mesh function $u$, the following relation holds at point $O\in \gamma_+$:
\begin{align}\label{eqn:average}
4u_O = u_{A}+u_{B}+u_{C}+u_{D}-\Delta_h u_O,
\end{align}
Points $A,B,C,D$ are its four connected points. At least one of them will belong to $M^+$ by definition of point set $\gamma$.


\paragraph{Case 0} Let's assume $x_A$ falls on the boundary, then $\alpha_A=0$ and $u(x_A)=u_A$, so
\begin{align}
4u_O = u(x_A)+u_C+u_D+u_E-\Delta_h u_O.
\end{align}
Not all points of $C,D,E$ will fall into $\gamma_+$. By the triangle inequality, we have
\begin{align}
4|u_{O}| \leq |u(x_A)|+|u_C|+|u_D|+|u_{E}| +|\Delta_h u_O|
\end{align}
and one of the following holds
\begin{subequations}
\begin{align}
4|u_{O}| &\leq |u|_{\Gamma_h} + 3|u|_\omega+|\Delta_h u|_{\gamma_+},\\
4|u_{O}| &\leq |u|_{\Gamma_h} + |u|_{\gamma_+} + 2|u|_\omega+|\Delta_h u|_{\gamma_+},\\
4|u_{O}| &\leq |u|_{\Gamma_h} + 2|u|_{\gamma_+} + |u|_\omega+|\Delta_h u|_{\gamma_+},\\
4|u_{O}| &\leq |u|_{\Gamma_h} + 3|u|_{\gamma_+}+|\Delta_h u|_{\gamma_+}.
\end{align}
\end{subequations}
In any case, the coefficients of $u_C,u_D,u_E$ satisfy
\begin{align}
4>1+1+1.
\end{align}

\begin{figure}[htbp]
    \centering
    \begin{subfigure}{0.3\textwidth}
        \centering
        \includegraphics[width=\textwidth]{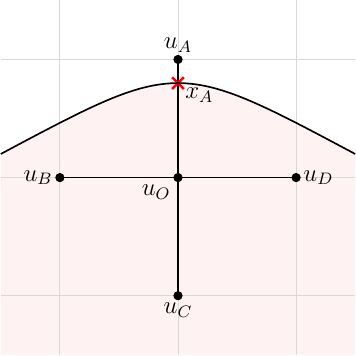}
        \caption{}
        \label{fig:one}
    \end{subfigure}
    ~
    \begin{subfigure}{0.3\textwidth}
        \centering
        \includegraphics[width=\textwidth]{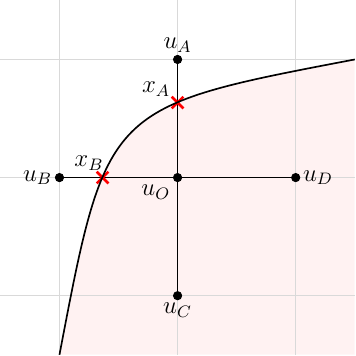}
        \caption{}
        \label{fig:two}
    \end{subfigure}
    ~
    \begin{subfigure}{0.3\textwidth}
        \centering
        \includegraphics[width=\textwidth]{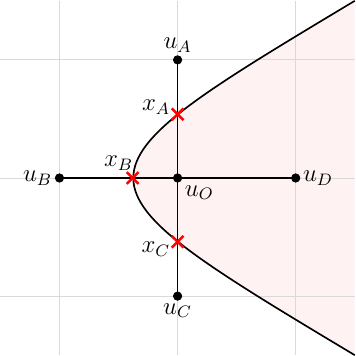}
        \caption{}
        \label{fig:three}
    \end{subfigure}
    \caption{Three cases of stencil points near the boundary (shaded area denotes the interior).}\label{fig:stencils}
\end{figure}

\paragraph{Case 1} 
Assume $A$ is the only exterior point as in Figure~\ref{fig:one}, the corresponding boundary point is $x_A$, and the rest three neighbor points $B,C,D$ are all inside the domain. The interpolation gives
\begin{align}
u(x_A) = \alpha_A u_A + (1-\alpha_A) u_O \Rightarrow u_A = \frac{1}{\alpha_A}u(x_A)+\left(1-\frac{1}{\alpha_A}\right)u_O.
\end{align}

Combining with the exact difference equation \eqref{eqn:average}, we have
\begin{align}
4u_O = \frac{1}{\alpha_A}u(x_A)+\left(1-\frac{1}{\alpha_A}\right)u_O + u_B+u_C+u_D-\Delta_h u_O
\end{align}
which can be simplified to 
\begin{align}
\left(3+\frac{1}{\alpha_A}\right)u_O = u_B+u_C+u_D+\frac{1}{\alpha_A}u(x_A),
\end{align}
This leads to one of the following:
\begin{subequations}
\begin{align}
\left(3+\frac{1}{\alpha_A}\right)|u_{O}| &\leq \frac{1}{\alpha_A}|u|_{\Gamma_h} + 3|u|_\omega+|\Delta_h u|_{\gamma_+},\\
\left(3+\frac{1}{\alpha_A}\right)|u_{O}| &\leq \frac{1}{\alpha_A}|u|_{\Gamma_h} + |u|_{\gamma_+} + 2|u|_\omega+|\Delta_h u|_{\gamma_+},\\
\left(3+\frac{1}{\alpha_A}\right)|u_{O}| &\leq \frac{1}{\alpha_A}|u|_{\Gamma_h} + 2|u|_{\gamma_+} + |u|_\omega+|\Delta_h u|_{\gamma_+},\\
\left(3+\frac{1}{\alpha_A}\right)|u_{O}| &\leq \frac{1}{\alpha_A}|u|_{\Gamma_h} + 3|u|_{\gamma_+}+|\Delta_h u|_{\gamma_+}.
\end{align}
\end{subequations}
Again, the coefficients of $u_O$ and its connected neighbors $u_B, u_C, u_D$ satisfy
\begin{align}
3+\frac{1}{\alpha_A}> 1+1+1.
\end{align}

\paragraph{Case 2} For two exterior points in Figure~\ref{fig:two}, we have linear interpolation
\begin{subequations}
\begin{align}
u_A &= \frac{1}{\alpha_A}u(x_A)+\left(1-\frac{1}{\alpha_A}\right)u_O,\\
u_B &= \frac{1}{\alpha_B}u(x_B)+\left(1-\frac{1}{\alpha_B}\right)u_O.
\end{align}
\end{subequations}
Plugging $u_A$ and $u_B$ into \eqref{eqn:average}, we obtain a similar expression for $u_O$:
\begin{align}
\left(2+\frac{1}{\alpha_A}+\frac{1}{\alpha_B}\right)u_O =u_C+u_D+\frac{1}{\alpha_A}u(x_A)+\frac{1}{\alpha_B}u(x_B)-\Delta_h u_O,
\end{align}
which again leads to one of the following
\begin{subequations}
\begin{align}
\left(2+\frac{1}{\alpha_A}+\frac{1}{\alpha_B}\right)|u_O| &\leq \left(\frac{1}{\alpha_A}+\frac{1}{\alpha_B}\right)|u|_{\Gamma_h} + 2|u|_{\omega}+|\Delta_h u|_{\gamma_+},\\
\left(2+\frac{1}{\alpha_A}+\frac{1}{\alpha_B}\right)|u_O| &\leq \left(\frac{1}{\alpha_A}+\frac{1}{\alpha_B}\right)|u|_{\Gamma_h} + |u|_{\gamma_+}+|u|_{\omega}+|\Delta_h u|_{\gamma_+},\\
\left(2+\frac{1}{\alpha_A}+\frac{1}{\alpha_B}\right)|u_O| &\leq \left(\frac{1}{\alpha_A}+\frac{1}{\alpha_B}\right)|u|_{\Gamma_h} + 2|u|_{\gamma_+}+|\Delta_h u|_{\gamma_+}.
\end{align}
\end{subequations}
The coefficients of $u_O$ and $u_C,u_D$ satisfy
\begin{align}
2+\frac{1}{\alpha_A}+\frac{1}{\alpha_B} > 1+1.
\end{align}

\paragraph{Case 3} Consider Figure~\ref{fig:three},
\begin{subequations}
\begin{align}
u_A &= \frac{1}{\alpha_A}u(x_A)+\left(1-\frac{1}{\alpha_A}\right)u_O,\\
u_B &= \frac{1}{\alpha_B}u(x_B)+\left(1-\frac{1}{\alpha_B}\right)u_O,\\
u_C &= \frac{1}{\alpha_B}u(x_B)+\left(1-\frac{1}{\alpha_B}\right)u_O.
\end{align}
\end{subequations}
Combining these interpolations with \eqref{eqn:average}, we obtain another form of the interpolation operator
\begin{align}
\left(1+\frac{1}{\alpha_A}+\frac{1}{\alpha_B}+\frac{1}{\alpha_C}\right)u_O=u_D+\frac{1}{\alpha_A}u(x_A)+\frac{1}{\alpha_B}u(x_B)+\frac{1}{\alpha_C}u(x_C)-\Delta_h u_O.
\end{align}
The inequality becomes
\begin{subequations}
\begin{align}
\left(1+\frac{1}{\alpha_A}+\frac{1}{\alpha_B}+\frac{1}{\alpha_C}\right)|u_O| &\leq \left(\frac{1}{\alpha_A}+\frac{1}{\alpha_B}+\frac{1}{\alpha_C}\right)|u|_{\Gamma_h} + |u|_{\omega}+|\Delta_h u|_{\gamma_+},\\
\left(1+\frac{1}{\alpha_A}+\frac{1}{\alpha_B}+\frac{1}{\alpha_C}\right)|u_O| &\leq \left(\frac{1}{\alpha_A}+\frac{1}{\alpha_B}+\frac{1}{\alpha_C}\right)|u|_{\Gamma_h} + |u|_{\gamma_+}+|\Delta_h u|_{\gamma_+},
\end{align}
\end{subequations}
and the coefficients of $u_O$ and $u_D$ satisfy
\begin{align}
1+\frac{1}{\alpha_A}+\frac{1}{\alpha_B}+\frac{1}{\alpha_C} > 1
\end{align}

\begin{figure}[htbp]
    \centering
    \begin{subfigure}{0.3\textwidth}
        \centering
        \includegraphics[width=\textwidth]{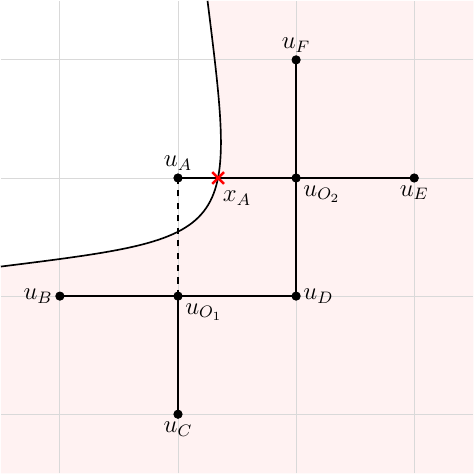}
        \caption{}
        \label{fig:1a}
    \end{subfigure}
    ~
    \begin{subfigure}{0.3\textwidth}
        \centering
        \includegraphics[width=\textwidth]{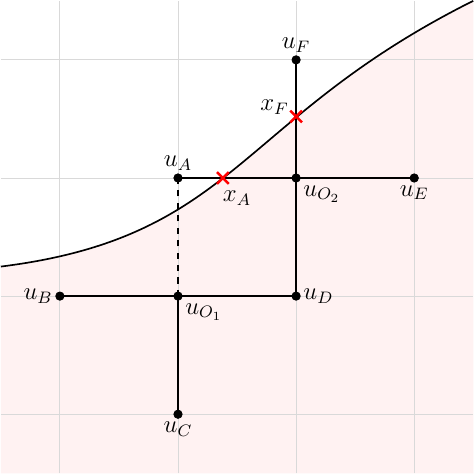}
        \caption{}
        \label{fig:1b}
    \end{subfigure}
    ~
    \begin{subfigure}{0.3\textwidth}
        \centering
        \includegraphics[width=\textwidth]{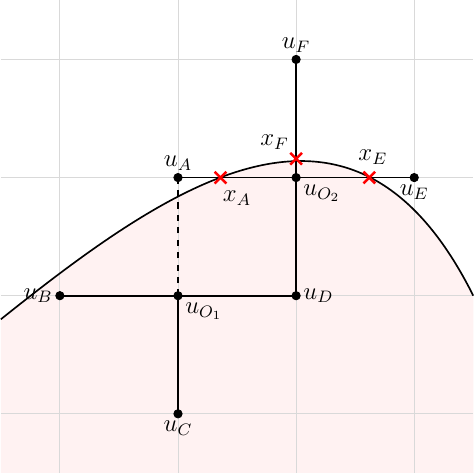}
        \caption{}
        \label{fig:1c}
    \end{subfigure}
    \caption{Three cases of stencil points near the boundary (shaded area denotes the interior).}\label{fig:stencil1}
\end{figure}

\paragraph{Case 1+X} In this case we explore the combinations of case 1 with cases 1, 2 and 3 as shown in Figure~\ref{fig:stencil1}.

First of all, at $O_1$, we have
\begin{align}
4u_{O_1} = \frac{1}{\alpha_A}u(x_A)+(1-\frac{1}{\alpha_A})u_{O_2}+u_B+u_C+u_D-\Delta_h u_{O_1}
\end{align}
and $u_{O_2}$ in the case of 1+1 (Figure~\ref{fig:1a}), 1+2 (Figure~\ref{fig:1b}), 1+3 (Figure~\ref{fig:1c}) satisfies, respectively
\begin{subequations}
\begin{align}
u_{O_2} &= \frac{1}{3+\frac{1}{\alpha_A}}\left(\frac{1}{\alpha_A}u(x_A)+u_D+u_E+u_F-\Delta_h u_{O_2}\right)\\
u_{O_2} &= \frac{1}{2+\frac{1}{\alpha_A}+\frac{1}{\alpha_F}}\left(\frac{1}{\alpha_A}u(x_A)+\frac{1}{\alpha_F}u(x_F)+u_D+u_E-\Delta_h u_{O_2}\right)\\
u_{O_2} &= \frac{1}{1+\frac{1}{\alpha_A}+\frac{1}{\alpha_F}+\frac{1}{\alpha_E}}\left(\frac{1}{\alpha_A}u(x_A)+\frac{1}{\alpha_F}u(x_F)+\frac{1}{\alpha_E}u(x_E)+u_D-\Delta_h u_{O_2}\right)
\end{align}
\end{subequations}
In all cases, after substitution of $u_{O_2}$ into the equation at $O_1$, the coefficients of the terms from $O_2$ are less than 1:
\begin{align}
\left|\frac{1-\frac{1}{\alpha_A}}{3+\frac{1}{\alpha_A}}\right|<1,\quad\left|\frac{1-\frac{1}{\alpha_A}}{2+\frac{1}{\alpha_A}+\frac{1}{\alpha_F}}\right|<1,\quad\left|\frac{1-\frac{1}{\alpha_A}}{1+\frac{1}{\alpha_A}+\frac{1}{\alpha_F}+\frac{1}{\alpha_E}}\right|<1
\end{align}

\begin{itemize}
\item In the case of 1+1, we have
\begin{align}
4u_{O_1} = \left(1+\frac{1-\frac{1}{\alpha_A}}{3+\frac{1}{\alpha_A}}\right)\frac{1}{\alpha_A}u(x_A)+u_B+u_C+\left(1+\frac{1-\frac{1}{\alpha_A}}{3+\frac{1}{\alpha_A}}\right)u_D\nonumber\\+\frac{1-\frac{1}{\alpha_A}}{3+\frac{1}{\alpha_A}}u_E+\frac{1-\frac{1}{\alpha_A}}{3+\frac{1}{\alpha_A}}u_F-\Delta_h u_{O_1}-\frac{1-\frac{1}{\alpha_A}}{3+\frac{1}{\alpha_A}}\Delta_h u_{O_2}
\end{align}

Note the absolute values of the coefficients satisfy
\begin{align}
4>1+1+1+\frac{1-\frac{1}{\alpha_A}}{3+\frac{1}{\alpha_A}}+\frac{-1+\frac{1}{\alpha_A}}{3+\frac{1}{\alpha_A}}+\frac{-1+\frac{1}{\alpha_A}}{3+\frac{1}{\alpha_A}}
\end{align}
The shared node $D$ cancels one contribution from $E$ or $F$.


\item In the case of 1+2,  
\begin{align}
4u_{O_1} = \left(1+\frac{1-\frac{1}{\alpha_A}}{2+\frac{1}{\alpha_A}+\frac{1}{\alpha_F}}\right)\frac{1}{\alpha_A}u(x_A)+\frac{1-\frac{1}{\alpha_A}}{2+\frac{1}{\alpha_A}+\frac{1}{\alpha_F}}\frac{1}{\alpha_F}u(x_F)+u_B+u_C\nonumber\\+\left(1+\frac{1-\frac{1}{\alpha_A}}{2+\frac{1}{\alpha_A}+\frac{1}{\alpha_F}}\right)u_D+\frac{1-\frac{1}{\alpha_A}}{2+\frac{1}{\alpha_A}+\frac{1}{\alpha_F}}u_E-\Delta_h u_{O_1}-\frac{1-\frac{1}{\alpha_A}}{2+\frac{1}{\alpha_A}+\frac{1}{\alpha_F}}\Delta_h u_{O_2}
\end{align}
The sum of the absolute values of coefficients of $u_B$, $u_C$, $u_D$ and $u_E$ also do not exceed 4
\begin{align}
4> 1+1+1+\frac{1-\frac{1}{\alpha_A}}{2+\frac{1}{\alpha_A}+\frac{1}{\alpha_F}}+\frac{-1+\frac{1}{\alpha_A}}{2+\frac{1}{\alpha_A}+\frac{1}{\alpha_F}}=3
\end{align}
thanks to the shared node $D$ again.

\item It is even better in the case of 1+3, since only points $B, C, D$ are involved.
\begin{align}
4u_{O_1} = \frac{1}{\alpha_A}u(x_A)+\frac{1-\frac{1}{\alpha_A}}{1+\frac{1}{\alpha_A}+\frac{1}{\alpha_F}+\frac{1}{\alpha_E}}\left(\frac{1}{\alpha_A}u(x_A)+\frac{1}{\alpha_F}u(x_F)+\frac{1}{\alpha_E}u(x_E)+u_D\right)\nonumber\\+u_B+u_C+u_D-\Delta_h u_{O_1}-\frac{1-\frac{1}{\alpha_A}}{1+\frac{1}{\alpha_A}+\frac{1}{\alpha_F}+\frac{1}{\alpha_E}}\Delta_h u_{O_2}
\end{align}
The sum of the absolute values of the coefficients of $u_B,u_C,u_D$ are also less than 4,
\begin{align}
4>3>1+1+1+\frac{1-\frac{1}{\alpha_A}}{1+\frac{1}{\alpha_A}+\frac{1}{\alpha_F}+\frac{1}{\alpha_E}}
\end{align}

\end{itemize}


\begin{figure}[htbp]
    \centering
    \begin{subfigure}{0.3\textwidth}
        \centering
        \includegraphics[width=\textwidth]{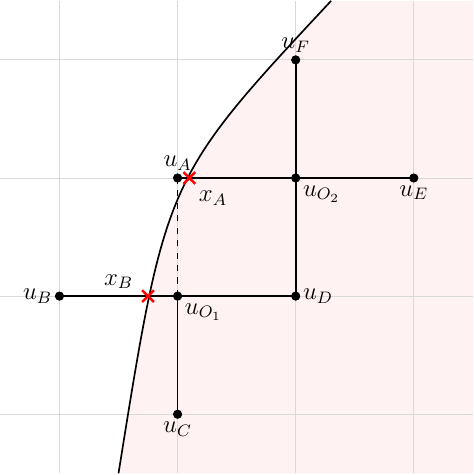}
        \caption{}
        \label{fig:2a}
    \end{subfigure}
    ~
    \begin{subfigure}{0.3\textwidth}
        \centering
        \includegraphics[width=\textwidth]{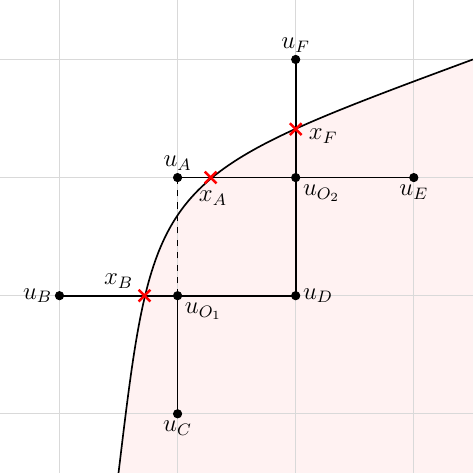}
        \caption{}
        \label{fig:2b}
    \end{subfigure}
    ~
    \begin{subfigure}{0.3\textwidth}
        \centering
        \includegraphics[width=\textwidth]{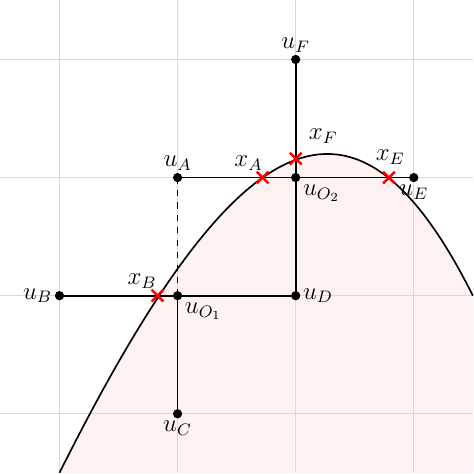}
        \caption{}
        \label{fig:2c}
    \end{subfigure}
    \caption{More cases of stencil points near the boundary (shaded area denotes the interior).}\label{fig:stencil2}
\end{figure}

\paragraph{Case 2+X} 
In Figure~\ref{fig:2a}, point $A$ chooses $O_2$ not $O_1$ to find the intersection point $x_A$. Now, the interpolation and
and the Laplace difference equation \eqref{eqn:average} give
\begin{align}
\left(3+\frac{1}{\alpha_B}\right)u_{O_1}=u_C+u_D+\frac{1}{\alpha_B}u(x_B)+\frac{1}{\alpha_A}u(x_A)+\left(1-\frac{1}{\alpha_A}\right)u_{O_2}-\Delta_h u_{O_1}
\end{align}
In the case of 2+1 (Figure~\ref{fig:2a}), 2+2 (Figure~\ref{fig:2b}) and 2+3 (Figure~\ref{fig:2c}), we have
\begin{subequations}
\begin{align}
\left(3+\frac{1}{\alpha_A}\right)u_{O_2} &= u_D+u_E+u_F+\frac{1}{\alpha_A}u(x_A)-\Delta_h u_{O_2}\\
\left(2+\frac{1}{\alpha_A}+\frac{1}{\alpha_F}\right)u_{O_2} &= u_D+u_E+\frac{1}{\alpha_A}u(x_A)+\frac{1}{\alpha_F}u(x_F)-\Delta_h u_{O_2}\\
\left(1+\frac{1}{\alpha_A}+\frac{1}{\alpha_E}+\frac{1}{\alpha_F}\right)u_{O_2} &= u_D+\frac{1}{\alpha_A}u(x_A)+\frac{1}{\alpha_E}u(x_E)+\frac{1}{\alpha_F}u(x_F)-\Delta_h u_{O_2}
\end{align}
\end{subequations}

We will not list all the inequalities, but only check the absolute values of the coefficients of the center point $O_1$ and points of $B,C,D,E,F,\dots$.
\begin{itemize}
    \item In the case of 2+1, the coefficients of $u_C,u_D,u_E,u_F$  satisfy
\begin{align}
3+\frac{1}{\alpha_B} > 3> 1+1+\frac{1-\frac{1}{\alpha_A}}{3+\frac{1}{\alpha_A}}+\frac{-1+\frac{1}{\alpha_A}}{3+\frac{1}{\alpha_A}}+\frac{-1+\frac{1}{\alpha_A}}{3+\frac{1}{\alpha_A}}
\end{align}

\item In the case of 2+2, the coefficients of $u_C, u_D, u_E$ is
\begin{align}
3+\frac{1}{\alpha_B}> 2=1+1+\frac{1-\frac{1}{\alpha_A}}{2+\frac{1}{\alpha_A}+\frac{1}{\alpha_F}}+\frac{-1+\frac{1}{\alpha_A}}{2+\frac{1}{\alpha_A}+\frac{1}{\alpha_F}}
\end{align}

\item In the case of 2+3, the coefficients of $D$ and $E$ satisfy
\begin{align}
3+\frac{1}{\alpha_B}>2>1+1+\frac{1-\frac{1}{\alpha_A}}{1+\frac{1}{\alpha_A}+\frac{1}{\alpha_F}+\frac{1}{\alpha_F}}
\end{align}

\end{itemize}
Thus, all cases of 2+X gives a dominant coefficient for $O_1$.

\paragraph{Case 3+X} First at $O_1$ we have
\begin{align}
\left(2+\frac{1}{\alpha_A}+\frac{1}{\alpha_B}\right)u_{O_1} = \frac{1}{\alpha_A}u(x_A)+\left(1-\frac{1}{\alpha_A}\right)u_{O_2}+\frac{1}{\alpha_B}u(x_B)+\frac{1}{\alpha_C}u(x_C)+u_D-\Delta_h u_{O_1}
\end{align}
Then in the case of 3+1 (Figure~\ref{fig:3a}), 3+2 (Figure~\ref{fig:3b}), 3+3 (Figure~\ref{fig:3c}), we have at $O_2$
\begin{subequations}
\begin{align}
\left(3+\frac{1}{\alpha_A}\right)u_{O_2}&=\frac{1}{\alpha_A}u(x_A)+u_F+u_E+u_D-\Delta_h u_{O_2}\\
\left(2+\frac{1}{\alpha_A}+\frac{1}{\alpha_F}\right)u_{O_2}&=\frac{1}{\alpha_A}u(x_A)+\frac{1}{\alpha_F}u(x_F)+u_E+u_D-\Delta_h u_{O_2}\\
\left(1+\frac{1}{\alpha_A}+\frac{1}{\alpha_F}+\frac{1}{\alpha_E}\right)u_{O_2}&=\frac{1}{\alpha_A}u(x_A)+\frac{1}{\alpha_F}u(x_F)+\frac{1}{\alpha_E}u(x_E)+u_D-\Delta_h u_{O_2}
\end{align}
\end{subequations}
respectively.

Now we check the coefficients at the grid points again.
\begin{itemize}
    \item 3+1: coefficients for $u_{O_1}$ and $u_D, u_E, u_F$
    \begin{align}
    2+\frac{1}{\alpha_A}+\frac{1}{\alpha_B} > 2 >  1+\frac{1-\frac{1}{\alpha_A}}{3+\frac{1}{\alpha_A}}+\frac{-1+\frac{1}{\alpha_A}}{3+\frac{1}{\alpha_A}}+\frac{-1+\frac{1}{\alpha_A}}{3+\frac{1}{\alpha_A}}
    \end{align}
    \item 3+2: coefficients for $u_{O_1}$ and $u_D, u_E$
    \begin{align}
    2+\frac{1}{\alpha_A}+\frac{1}{\alpha_B} > 2 >  1+\frac{1-\frac{1}{\alpha_A}}{2+\frac{1}{\alpha_A}+\frac{1}{\alpha_F}}+\frac{-1+\frac{1}{\alpha_A}}{2+\frac{1}{\alpha_A}+\frac{1}{\alpha_F}}=1
    \end{align}
    \item 3+3: coefficients for $u_{O_1}$ and $u_D$
    \begin{align}
    2+\frac{1}{\alpha_A}+\frac{1}{\alpha_B} > 1 >  1+\frac{1-\frac{1}{\alpha_A}}{1+\frac{1}{\alpha_A}+\frac{1}{\alpha_F}+\frac{1}{\alpha_E}}
    \end{align}

\end{itemize}

\begin{figure}[htbp]
    \centering
    \begin{subfigure}{0.3\textwidth}
        \centering
        \includegraphics[width=\textwidth]{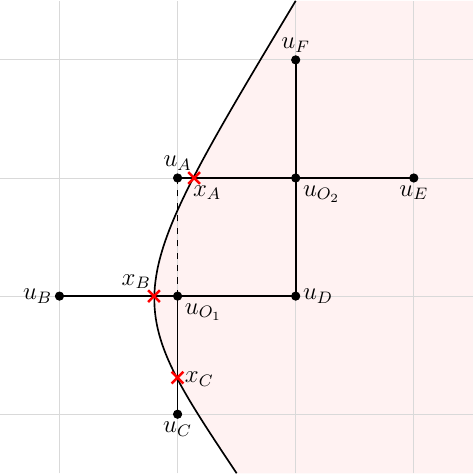}
        \caption{}
        \label{fig:3a}
    \end{subfigure}
    ~
    \begin{subfigure}{0.3\textwidth}
        \centering
        \includegraphics[width=\textwidth]{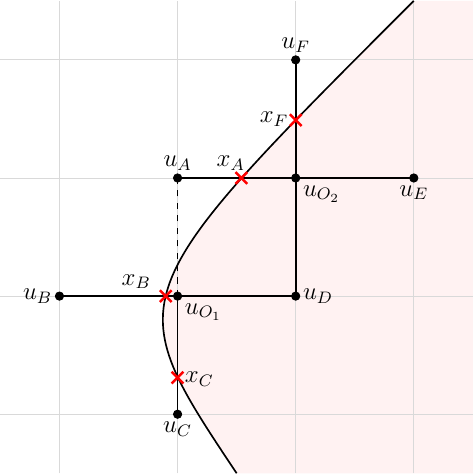}
        \caption{}
        \label{fig:3b}
    \end{subfigure}
    ~
    \begin{subfigure}{0.3\textwidth}
        \centering
        \includegraphics[width=\textwidth]{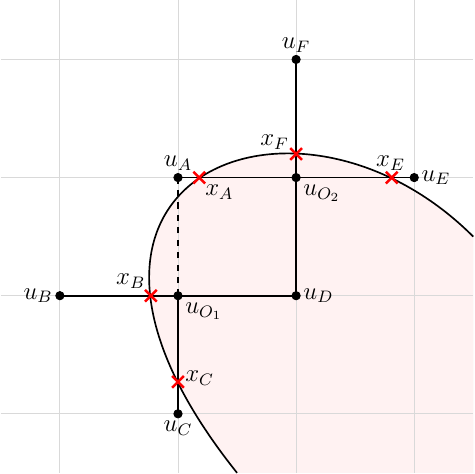}
        \caption{}
        \label{fig:3c}
    \end{subfigure}
    \caption{More cases of stencil points near the boundary (shaded area denotes the interior).}\label{fig:stencil3}
\end{figure}

\paragraph{Case $N$-Chain} There are numerous combinations of the cases discussed above and the number of coupled ($O_1,O_2,\dots,O_n$) points can be arbitrary. The end points of the chain can be any of cases 1, 2 and 3, and inside the chain, only cases 2 and 3 are possible. In Figure~\ref{fig:stencilN}, we show the case of a 3-Chain: 2+2+2.

\begin{figure}[htbp]
\centering
\includegraphics[width=0.35\textwidth]{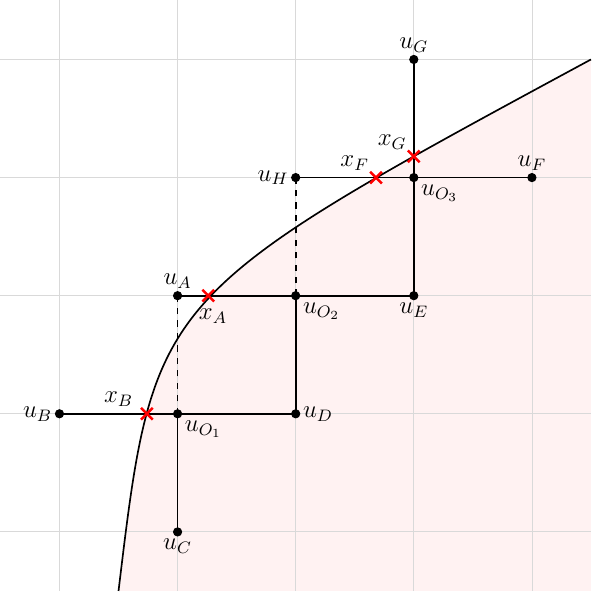}
\caption{$N$-Chain}\label{fig:stencilN}
\end{figure}

At each point $O_1, O_2, O_3$, we have
\begin{subequations}
\begin{align}
\left(3+\frac{1}{\alpha_B}\right)u_{O_1} &= \frac{1}{\alpha_B}u(x_B)+u_C+u_D+\frac{1}{\alpha_A}u(x_A)+(1-\frac{1}{\alpha_A})u_{O_2}-\Delta_h u_{O_1}\\
\left(3+\frac{1}{\alpha_A}\right)u_{O_2} &= \frac{1}{\alpha_A}+u_D + u_E + \frac{1}{\alpha_H}u(x_H)+(1-\frac{1}{\alpha_H})u_{O_3}-\Delta_h u_{O_2}\\
\left(2+\frac{1}{\alpha_A}+\frac{1}{\alpha_B}\right)u_{O_3} &= u_E+u_F+\frac{1}{\alpha_H}u(x_H)+\frac{1}{\alpha_G}u(x_G)-\Delta_h u_{O_3}.
\end{align}
\end{subequations}
We can substitute $u_{O_3}$ into the equation of $u_{O_2}$, and $u_{O_2}$ into the equation of $u_{O_1}$. Now we can check the coefficients of $u_{O_1}$ and interior nodes used $u_C,u_D,u_E,u_F$.
\begin{align}
3+\frac{1}{\alpha_B}>1+\left(1+\frac{1-\frac{1}{\alpha_A}}{3+\frac{1}{\alpha_A}}\right)+\left(\frac{-1+\frac{1}{\alpha_A}}{3+\frac{1}{\alpha_A}}\right)\left(1+\frac{1-\frac{1}{\alpha_H}}{2+\frac{1}{\alpha_H}+\frac{1}{\alpha_G}}\right)+\frac{1-\frac{1}{\alpha_A}}{3+\frac{1}{\alpha_A}}\frac{1-\frac{1}{\alpha_H}}{2+\frac{1}{\alpha_H}+\frac{1}{\alpha_G}}=2
\end{align}

The shared nodes $D$ and $E$ cancel some of the terms, making the coefficients of $u_{O_1}$ dominant. Other chains of more than two central points ($O_1,O_2,\dots,O_n$) follow the arguments above similarly. 

Thus, for any point $O\in\gamma_+$, we can conclude that
\begin{align}
|u_O| \leq  \frac{\alpha_\omega}{\alpha_O}|u|_{\omega}+\frac{\alpha_{\gamma_+}}{\alpha_O}|u|_{\gamma_+}+\frac{\alpha_{\Gamma_h}}{\alpha_O}|\tau|_{\Gamma_h}+\frac{\alpha_\Delta}{\alpha_O}|\Delta u|_{\gamma_+}
\end{align}
where $\alpha_O>\alpha_\omega + \alpha_{\gamma_+}$. Since point $O$ is arbitrary, we can conclude that
\begin{align}
|u|_{\gamma} \leq \max_{O\in\gamma_+}\left\{\frac{\alpha_\omega}{\alpha_O}|u|_{\omega}+\frac{\alpha_{\gamma_+}}{\alpha_O}|u|_{\gamma_+}+\frac{\alpha_{\Gamma_h}}{\alpha_O}|\tau|_{\Gamma_h}+\frac{\alpha_\Delta}{\alpha_O}|\Delta u|_{\gamma_+}\right\}
\end{align}
and we let the maximum case be
\begin{align}
|u|_{\gamma} \leq \beta_\omega|u|_{\omega}+\beta_{\gamma_+}|u|_{\gamma_+}+\beta_{\Gamma_h}|\tau|_{\Gamma_h}+\beta_\Delta|\Delta u|_{\gamma_+}
\end{align}
which leads to
\begin{align}
|u|_{\gamma} \leq \frac{\beta_\omega}{1-\beta_{\gamma_+}}|u|_{\omega}+\frac{\beta_{\Gamma_h}}{1-\beta_{\gamma_+}}|\tau|_{\Gamma_h}+\frac{\beta_{\Delta}}{1-\beta_{\gamma_+}}|\Delta u|_{\gamma_+}
\end{align}
Since $1>\beta_\omega+\beta_{\gamma_+}$, we can conclude that \eqref{eqn:ineq} holds and $\beta_1$ in \eqref{eqn:ineq} is less than 1.
\end{proof}

\begin{remark}
When small cut occurs, the $\alpha$'s will be close to 0, which will make $1/\alpha$ large and will not introduce instability.
\end{remark}

\begin{remark}
The interpolation operator uses only interior points and boundary points for interpolation at points in $\gamma_+$. Thus, no extrapolation is used, even though our discretization of the boundary condition looks like extrapolation.
\end{remark}

\begin{remark}
If we will focus on difference Laplace equation at $m\in \gamma_+$, i.e.
\begin{align}
Pu_{\gamma_+} = Qu_{w}+Ru_{\Gamma_h},
\end{align}
then
\begin{align}
||P^{-1}Q||_{\infty}\leq 1.
\end{align}
\end{remark}

Following \cite{larsson2003partial}, we give the interior estimate before proving the accuracy.
\begin{lemma}\label{lem:interior}
The interior estimate holds for any mesh function $u$.
\begin{align}
|u|_{M^+}\leq C_1|\tau|_{\Gamma_h}+C_2|\Delta u|_{M^+}
\end{align}
\end{lemma}

\begin{proof}
First note that $\Delta_hu=0$ at all points in $M^+$, by maximum principle, we have
\begin{align}
|u|_{M^+} \leq |u|_{\gamma_+},
\end{align}
if one regards $\gamma_+$ as the discrete grid boundary of set $M^+$.

Also,
\begin{align}
|u|_{\omega} \leq |u|_{M^+}, \quad |\Delta_h u|_{\gamma_+}\leq |\Delta u|_{M^+}
\end{align}
since $\omega,\gamma_+\subset M^+$.

By Lemma~\ref{lem:stability}, we reach
\begin{align}
|u|_{M^+}\leq |u|_{\gamma_+} \leq \beta_1|u|_\omega+\beta_2|\tau|_{\Gamma_h}+\beta_3|\Delta_h u|_{\gamma_+} \leq \beta_1|u|_{M^+}+\beta_2|\tau|_{\Gamma_h}+\beta_3|\Delta_h u|_{M^+}
\end{align}
Since $\beta_1<1$, we have
\begin{align}
(1-\beta_1)|u|_{M^+} \leq \beta_2|\tau|_{\Gamma_h}+\beta_3|\Delta_h u|_{M^+}\Rightarrow |u|_{M^+} \leq \frac{\beta_2}{1-\beta_1}|\tau|_{\Gamma_h}+\frac{\beta_3}{1-\beta_1}|\Delta_h u|_{M^+}
\end{align}
\end{proof}

\begin{theorem}
Let $u_{exact}$ be the solution of Laplace equation with boundary condition $u=g$ and $u$ be our construction using the boundary algebraic equation approach, then
\begin{align}
|u_{exact}-u|_{M^+} \leq C h^2
\end{align} 
\end{theorem}
\begin{proof}
The error $z_j = u_j-u_{exact}(x_j)$ satisfies
\begin{align}
\Delta_hz_j = 0-\Delta_h u_{exact}(x_j)=\Delta u_{exact}(x_j)-\Delta_hu_{exact}(x_j)\leq Ch^2
\end{align}
Applying Lemma~\ref{lem:interior}, we have
\begin{align}
|z|_{M^+}\leq C_1|\tau(z)|_{\Gamma_h}+C_2|\Delta_h z|_{M^+}\leq Ch^2,
\end{align}
since the interpolation is also second order accurate.
\end{proof}

\begin{remark}
The Neumann or Robin BC case can be argued similarly, as long as we can resolve the boundary conditions to second order accuracy, the convergence follows from the stability.
\end{remark}

\begin{remark}
Consider the 1D difference Laplace equation with Dirichlet BC at $x=a$ and $x=b$:
\begin{subequations}
\begin{align}
u_{i-1}-2u_i+u_{i+1} = 0, i=1,2,\dots,N\\
\alpha u_0 + (1-\alpha) u_1 = g(a)\\
\beta u_{N+1} + (1-\beta) u_N = g(b)
\end{align}
\end{subequations}
where $0<\alpha,\beta<1$.
Using the equation at $x_1$: $u_0-2u_1+u_2=0$ we have 
\begin{subequations}
\begin{align}
\alpha (2u_1-u_2) + (1-\alpha) u_1 = g(a) \\
 \Rightarrow u_1 = \frac{\alpha}{1+\alpha}u_2+\frac{1}{1+\alpha}g(a)
\end{align}
\end{subequations}
which means that the closet interior point to the boundary is a convex linear combination of its neighboring points, one exterior and one interior point, which defines exactly the interpolation operator $\ell_h$ used in Larsson and Thom\'{e}e's book \cite[page 48]{larsson2003partial} for unfitted finite difference. The stability and accuracy follow closely Lemma 4.5 and Theorem 4.3 \cite[page 48]{larsson2003partial}. Our 2D unfitted finite difference scheme and the stability argument bear similar argument as the 1D case and can be regarded as a generalization of those in \cite{larsson2003partial}.
\end{remark}
\section{Fast computation of homogeneous solutions via difference potentials}\label{sec:dpm}

Recall that once we find the single layer density $q_s$ or the double layer density $q_d$, we will be able to compute the value $u(m)$ via the convolutions
\begin{align}\label{eqn:convolutions}
u(m) = \sum_n S(m,n)q_s(n) = \sum_n D(m,n)q_d(n),
\end{align}
for any interior point $x_m\in M^+$. However, summing up for every point in $M^+$ would need kernel information of the source-target interaction and the summation cost would be in the scale of $\mathcal{O}(N^3)$.

To speed up the summation in bounded domains, we follow concepts in difference potentials method \cite{ryaben2012method} and only compute $u_\gamma$ for points $x_m\in\gamma$ using the convolutions \eqref{eqn:convolutions} with $\mathcal{O}(N^2)$ cost. 

Next, as in the difference potentials framework, we introduce the following auxiliary problem defined on a larger rectangular box $\Omega^0\supset\Omega$.
\begin{subequations}\label{eqn:aux_dp}
\begin{align}
[Aw](m) &= \left\{
\begin{array}{cc}
0, & m\in M^+,\\
Au, & m\in M^-,\\
\end{array}
\right.\\
w(m) &= 0,\quad m\in \partial\Omega^0.
\end{align}
\end{subequations}
where
\begin{align}\label{eqn:zero_extension}
u = \left\{
\begin{array}{cc}
u_\gamma(m), & m\in\gamma,\\
0, & \mbox{elsewhere}.\\
\end{array}
\right.
\end{align}
Here, the obtained $u_\gamma$ is zero-padded into $u$. Clearly, the solution $w$ satisfies the Laplace equation in $M^+$ and we show that the trace of $w$ on $\gamma$ is exactly $u_\gamma$. Note that the auxiliary problem \eqref{eqn:aux_dp} can be solved using FFT based fast solvers with complexity $\mathcal{O}(N^2\log(N))$.

\begin{theorem}\label{thm:gamma}
The solution $w$, also known as the difference potentials when restricted to $N^+$, to the system \eqref{eqn:aux_dp} solves the homogeneous difference equation $Aw=0$ in $M^+$ and coincides with $u_\gamma$ on point set $\gamma$, i.e.
\begin{align}
Tr_{\gamma}w = u_\gamma,
\end{align}
where $Tr_{\gamma}$ denotes the trace operator.
\end{theorem}

\begin{proof}
First of all, we define point set $M^+,M^-,N^+,N^-,\gamma$ inside $\Omega^0$ similarly, except the domain is bounded now. Then $M^0:=M^+\cup M^-$ and $N^0:=N^+\cup N^-$. Equation \eqref{eqn:aux_dp} admits a unique solution $w$ and $w$ satisfies the homogeneous difference equation $Aw=0$ in $M^+$ by construction.

Next consider the follow system defined on $N^+$ given $u_\gamma$
\begin{subequations}
\begin{align}
[Av_{N^+}](m) &=0,\quad m\in M^+,\\
v_{N^+}(m) &= u_\gamma(m), \quad m \in \gamma.
\end{align}
\end{subequations}
By construction, clearly we have $Tr_\gamma v_{N^+} = u_\gamma$.

Now extend $v_{N^+}$ to $v_{N^0}$ by 0 from $N^+$ to $N^0$, thus $v_{N^0}$ satisfies the homogeneous boundary condition on $\partial \Omega^0$.

Lastly consider the difference between $w$ and $v_{N^0}$, then we have
\begin{subequations}\label{eqn:aux_diff}
\begin{align}
[A(w-v_{N^0})](m) &= \left\{
\begin{array}{cc}
0, & m\in M^+,\\
A(u-v_{N^0}), & m\in M^-,\\
\end{array}
\right.\\
w(m) &= 0,\quad m\in \partial\Omega^0.
\end{align}
\end{subequations}
For $m\in N^-$, $u(m)=v_{N^0}(m)$, then \eqref{eqn:aux_diff} becomes
\begin{subequations}
\begin{align}
[A(w-v_{N^0})](m) &= 0,\quad m\in M^0,\\
w(m) &= 0,\quad m\in \partial\Omega^0.
\end{align}
\end{subequations}
which admits only the trivial solution, thus
\begin{align}
Tr_{\gamma} w = Tr_{\gamma}v_{N^0} = Tr_{\gamma}v_{N^+} = u_\gamma,
\end{align}
which concludes the proof.
\end{proof}

Based on Theorem~\ref{thm:gamma}, we can solve the difference equation~\eqref{eqn:aux_dp} once to retrieve all interior nodal values in $M^+$ if we know the boundary nodal values $u_\gamma$. This can also be interpreted as that discrete convolution with single layer kernel can be evaluated using forward and inverse Fourier transforms, which is equivalent to using fast Poisson solver based on fast Fourier Transform to solve the system \eqref{eqn:aux_dp} for speedup.

\begin{remark}
Theorem~\ref{thm:gamma} is a classic result in the difference potentials framework \cite{ryaben2012method}, that is, for homogeneous difference equations, the trace of the difference potentials is same as the density $u_\gamma$. Essentially, the difference potentials operator is a projection operator.
\end{remark}

\section{Nonhomogeneous equations}\label{sec:nonhomogeneous_equations}
For the following nonhomogeneous difference equation
\begin{align}\label{eqn:nonhomo}
[Au](m): = 4u(m)-u(m+e_1)-u(m-e_1)-u(m+e_2)-u(m-e_2) = h^2f(m),
\end{align}
we decompose $u(m)$ into a homogeneous solution $u^h(m)$ which satisfies the homogeneous difference equation~\eqref{eqn:homo} and a particular solution $u^p(m)$ that satisfies the nonhomogeneous difference equation.

The choice of the particular solution is arbitrary. One can, for example, use the convolution of forcing terms with the LGF
\begin{align}
u^p(m) = \sum_{n\in M^+} h^2f(n)K(m,n)
\end{align}
where the kernel $K$ is either the single layer kernel $S$ or the double layer kernel $D$. When $f$ is compactly supported, this approach is viable and fast summation techniques based fast multipole method such as \cite{gillman2014fast} would be favorable.



Another approach is to employ an auxiliary problem again. We solve for a particular solution $u^p(m)$ on mesh $\Omega^0_h$ with homogeneous boundary condition at the boundary of $\Omega^0$, i.e.
\begin{subequations}\label{eqn:ps}
\begin{align}
[Au^p](m) &= \left\{
\begin{array}{cc}
h^2f(m), & m\in M^+,\\
0, & m\in M^-,\\
\end{array}
\right.\\
u^p(m) &= 0,\quad m\in \partial\Omega^0.
\end{align}
\end{subequations}

\begin{remark}
The boundary condition in \eqref{eqn:ps} is arbitrary as long as it's well-posed. Here, we choose homogeneous boundary condition to allow for fast Poisson solver based on fast Fourier transform, so that we can solve for the particular solution rapidly.
\end{remark}

Once we obtain information about the particular solution $u^p(m)$, we can subtract off the contribution of the particular solution in the discretization of the boundary condition. For example, in the Dirichlet case,
\begin{align}
\sum_{x_{jk}\in\gamma} (u^p_{jk}+u^h_{jk})\phi_{jk}(x_i) = g(x_i).
\end{align}
Since $u^p$ is known, we have
\begin{align}
\sum_{x_{jk}\in\gamma} u^h_{jk}\phi_{jk}(x_i) = g(x_i)-\sum_{x_{jk}\in\gamma} u^p_{jk}\phi_{jk}(x_i).
\end{align}
Or in the matrix-vector form, we have
\begin{align}
\Phi u^h_{\gamma} = g(x_i)-\Phi u^p_{\gamma}.
\end{align}

The Neumann or Robin boundary condition can be treated similarly by subtracting off contribution in the boundary conditions from the particular solution.

Next, we follow the homogeneous case and solve for the homogeneous solution $u^h(m)$. Lastly, using superposition, we add the homogeneous solution to the particular solution, which gives the approximated solution $u(m) = u^h(m) + u^p(m)$ for points in $M^+$.
\section{Numerical Results}\label{sec:num}
Unless specified otherwise, we test with manufactured solution
\begin{align}
u(x,y) = \sin(x)\cos(y)
\end{align}
and the forcing functions, Dirichlet BC or Robin BC are computed using the manufactured solution. The computational domain is selected to be $[-1-\ell,1+\ell]\times [-1-\ell,1+\ell]$ with $\ell=0.15$ unless specified otherwise. The grid size $h$ then is given by $h=(2+2\ell)/N$. The convergence is computed with errors in max norm.

\subsection{Dirichlet BC in ellipses}
We first test the developed method with Dirichlet BC and ellipses of different aspect ratios $\alpha$:
\begin{align}
x^2+\alpha^2y^2=1.
\end{align}

\begin{figure}[htbp]
    \centering
    \begin{subfigure}{0.45\textwidth}
        \centering
        \includegraphics[width=\textwidth, trim = 2cm 7cm 2cm 6.5cm]{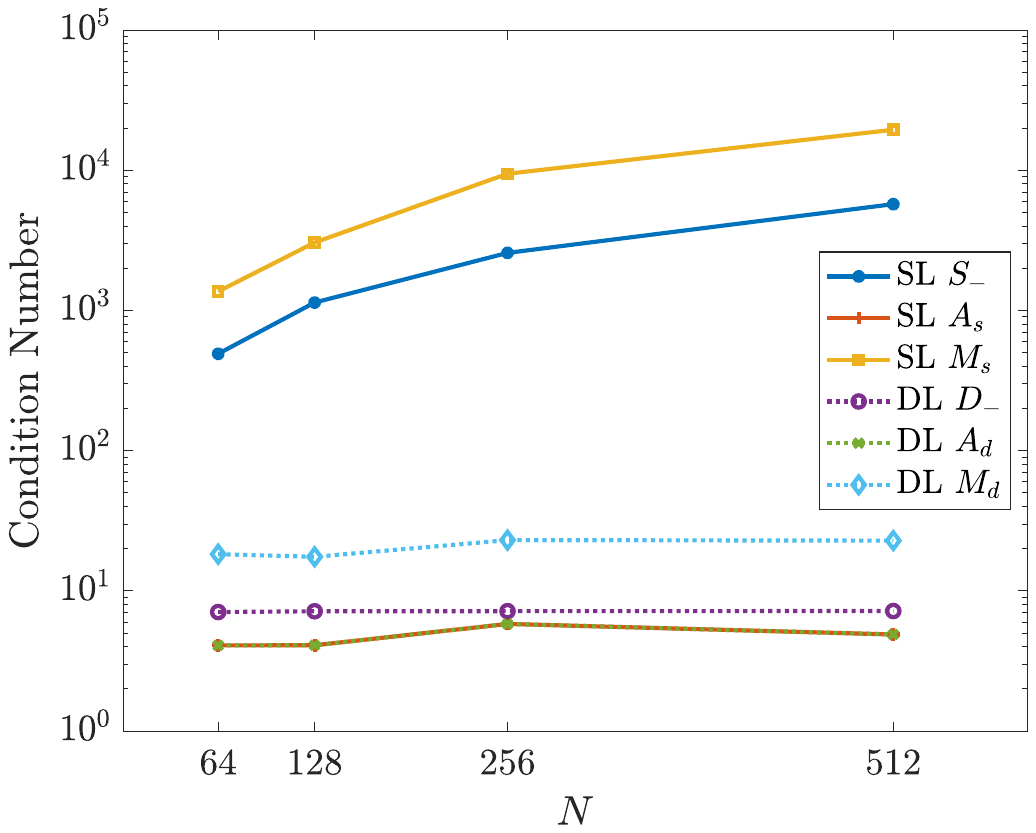}
        \caption{$\alpha=2$}
    \end{subfigure}
    ~
    \begin{subfigure}{0.45\textwidth}
        \centering
        \includegraphics[width=\textwidth, trim = 2cm 7cm 2cm 6.5cm]{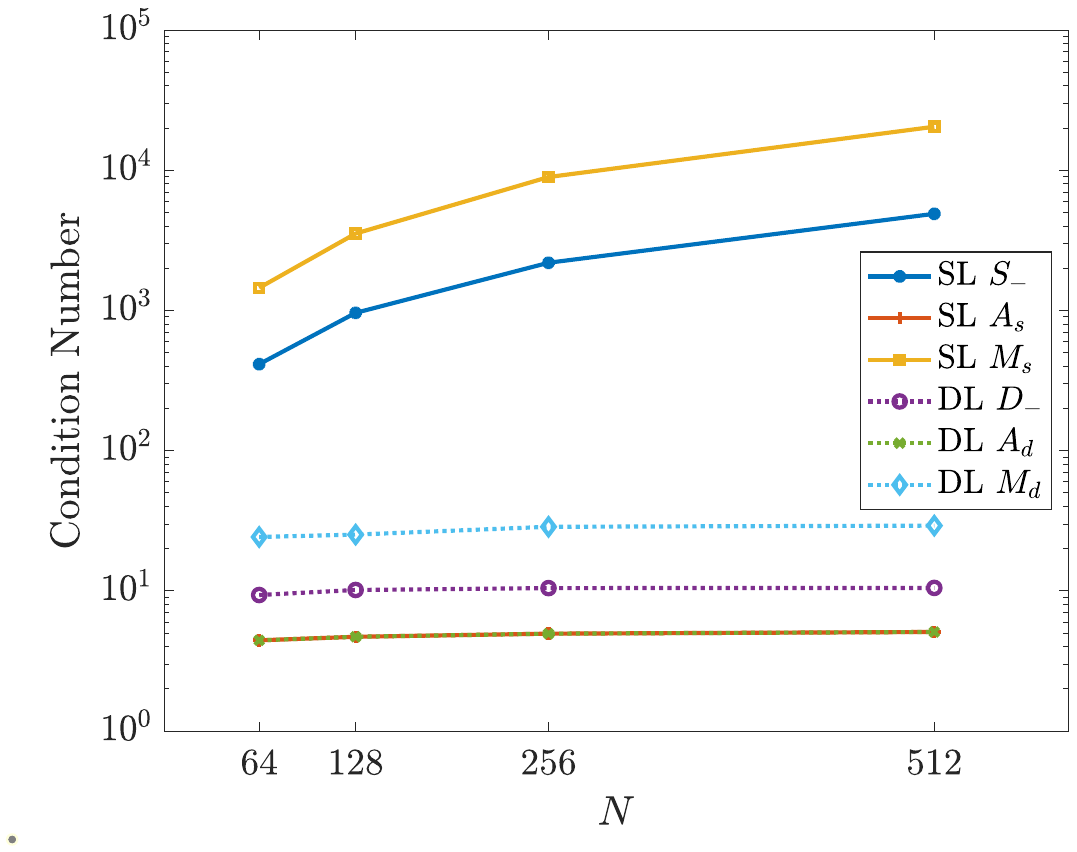}
        \caption{$\alpha=8$}
    \end{subfigure}
    \caption{Condition number (Dirichlet BC) for ellipses of different aspect ratios}\label{fig:ellipse_cond}
\end{figure}

The matrices that we compute the condition numbers are defined as
\begin{align}
A_s&=\Phi_+S_+S^{-1}_{-} +\Phi_-,\\
M_s&=\Phi_+S_+ +\Phi_-S_{-},\\
A_d&=\Phi_+D_+D^{-1}_{-} +\Phi_-,\\
M_d&=\Phi_+D_+ +\Phi_-D_{-}.
\end{align}
In Fig~\ref{fig:ellipse_cond}, it is observed that the condition numbers are independent of the aspect ratios of the geometry where we only show two cases of $\alpha=2, 8$. In particular, the single layer formulations $S_-$ and $M_s$ demonstrate mild and logarithmic growth under mesh refinement, while $A_s$ (which is preconditioned) and double layer formulations $D_-, A_d, M_d$ are all uniformly bounded under mesh refinement, which is desirable for iterative solvers such as GMRES. It is also interesting to note that the preconditioned versions $A_s$ and $A_d$ share identical condition number in all ellipses of different aspect ratios.

\begin{figure}[htbp]
\centering
    \begin{subfigure}{0.45\textwidth}
        \centering
        \includegraphics[width=\textwidth, trim = 2cm 7cm 2cm 6.5cm]{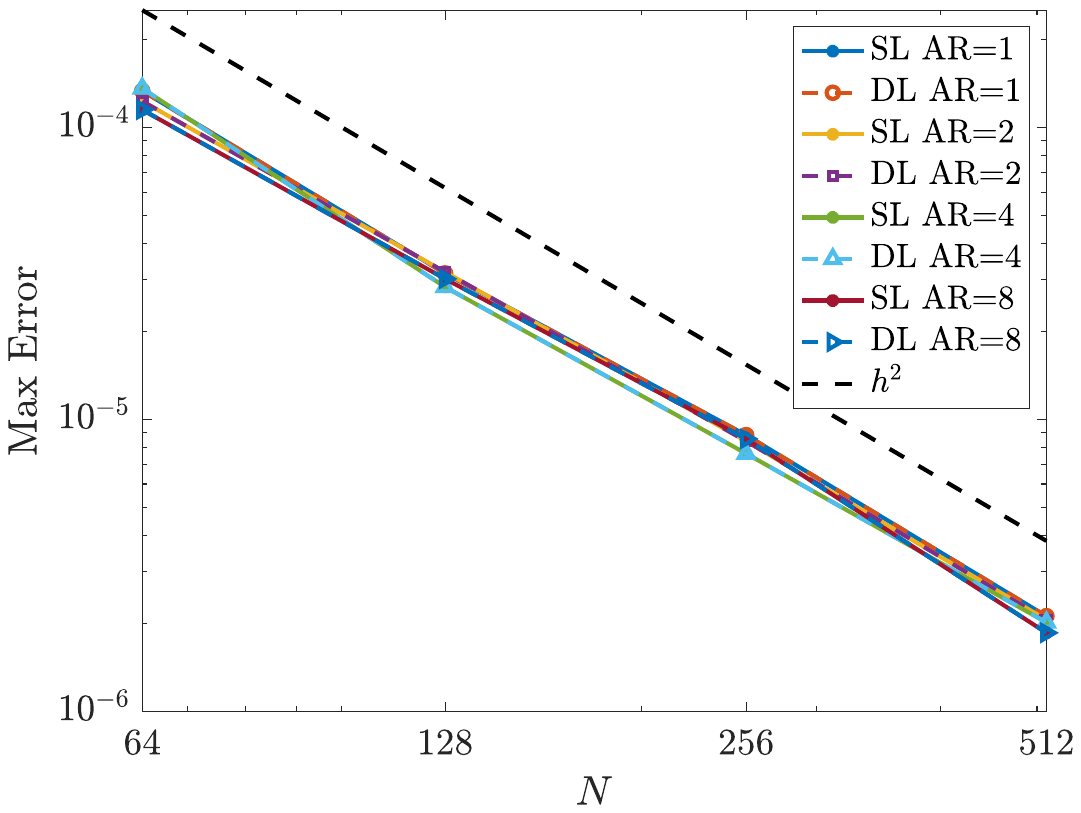}
        \caption{Convergence}\label{fig:ellipse_d_conv}
    \end{subfigure}
    ~
    \begin{subfigure}{0.45\textwidth}
        \centering
        \includegraphics[width=\textwidth, trim = 2cm 7cm 2cm 6.5cm]{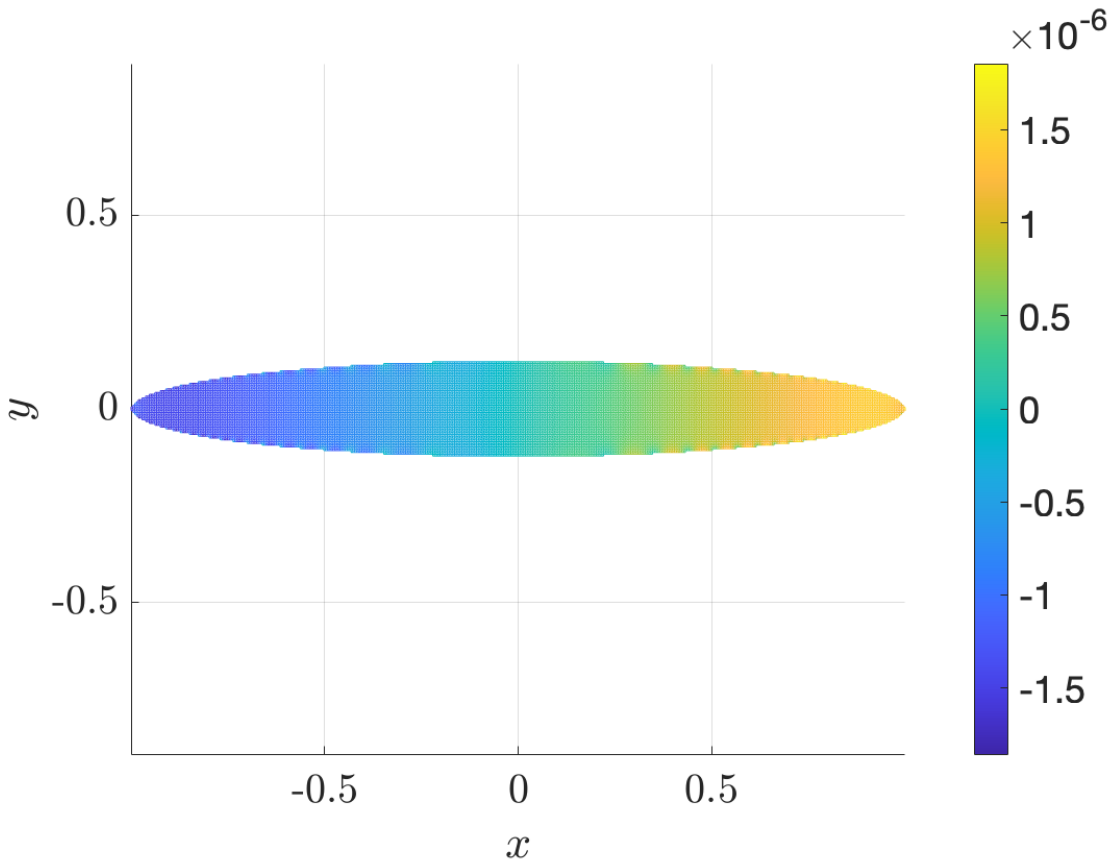}
        \caption{Error pattern for $\alpha=8$ ($512\times 512$)}\label{fig:ellipse_d_error}
    \end{subfigure}
    \caption{Convergence (Dirichlet BC) for ellipses of different aspect ratios and errors for aspect ratio $\alpha=8$}\label{fig:ellipse_d}
\end{figure}

In Fig~\ref{fig:ellipse_d}, we observe second order convergence for all ellipses with different aspect ratios (1,2,4,8). The magnitudes of the max errors do not depend on the aspect ratios of the geometry, indicating robustness of the developed method for different shapes of geometry. Neither do the max errors depend on the single or double layer formulations. The errors in Fig~\ref{fig:ellipse_d_error} show that the errors are smooth and the errors has large absolute values where the solution itself is large.

\subsection{Dirichlet BC in diamonds}
Next, we test the Dirichlet BC in diamonds of different aspect ratios. Overall, the numerical result do not differ much from the cases of ellipses of different aspect ratios.

The diamond shape is given by the zero level set of the following function
\begin{align}
\psi(x,y) = \left|\frac{x}{r_1}\right|+\left|\frac{y}{r_2}\right|-1
\end{align}
where $r_1=0.9$ and $r_2=0.5$. For such geometries with corners, boundary integral method requires special quadrature such as those developed in \cite{helsingfast2011,bremer2010efficient,helsing2013solving}, whereas our method based on discrete potentials treat smooth geometries and piecewise smooth geometries alike. This is due to the fact that we enforce boundary conditions at intersection points of the boundary and the mesh. 

\begin{figure}[htbp]
    \centering
    \begin{subfigure}{0.45\textwidth}
        \centering
        \includegraphics[width=\textwidth, trim = 2cm 7cm 2cm 6.5cm]{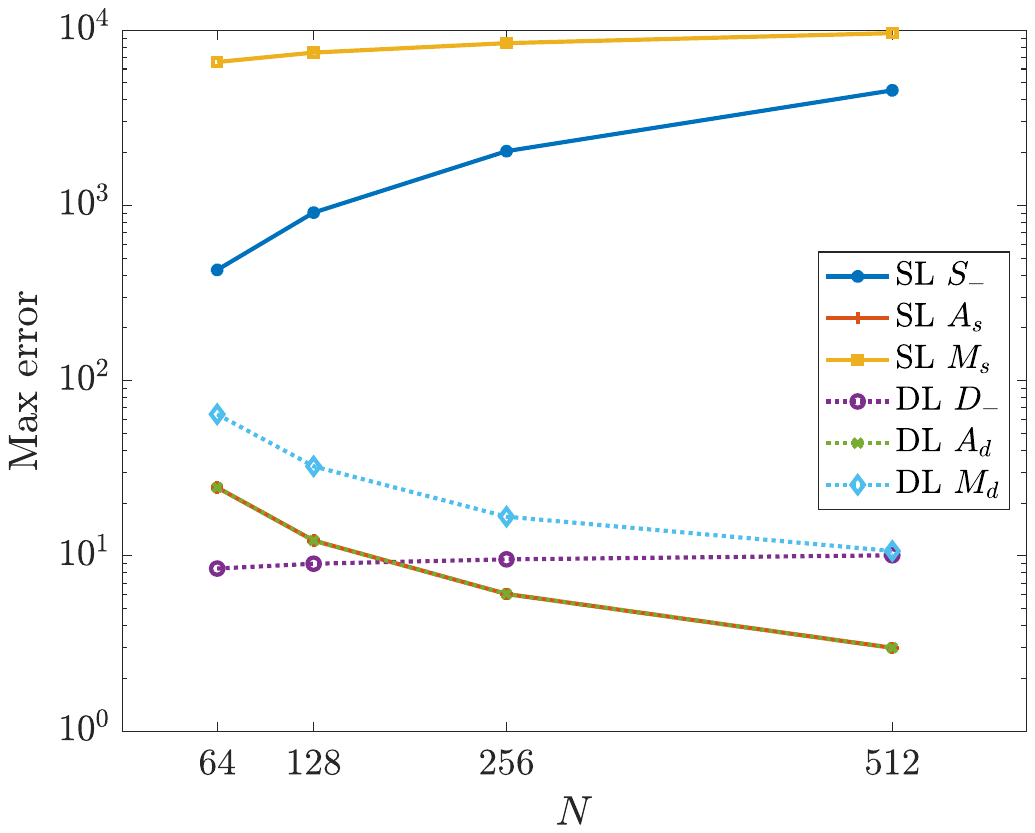}
        \caption{$\alpha=2$}
    \end{subfigure}
    ~
    \begin{subfigure}{0.45\textwidth}
        \centering
        \includegraphics[width=\textwidth, trim = 2cm 7cm 2cm 6.5cm]{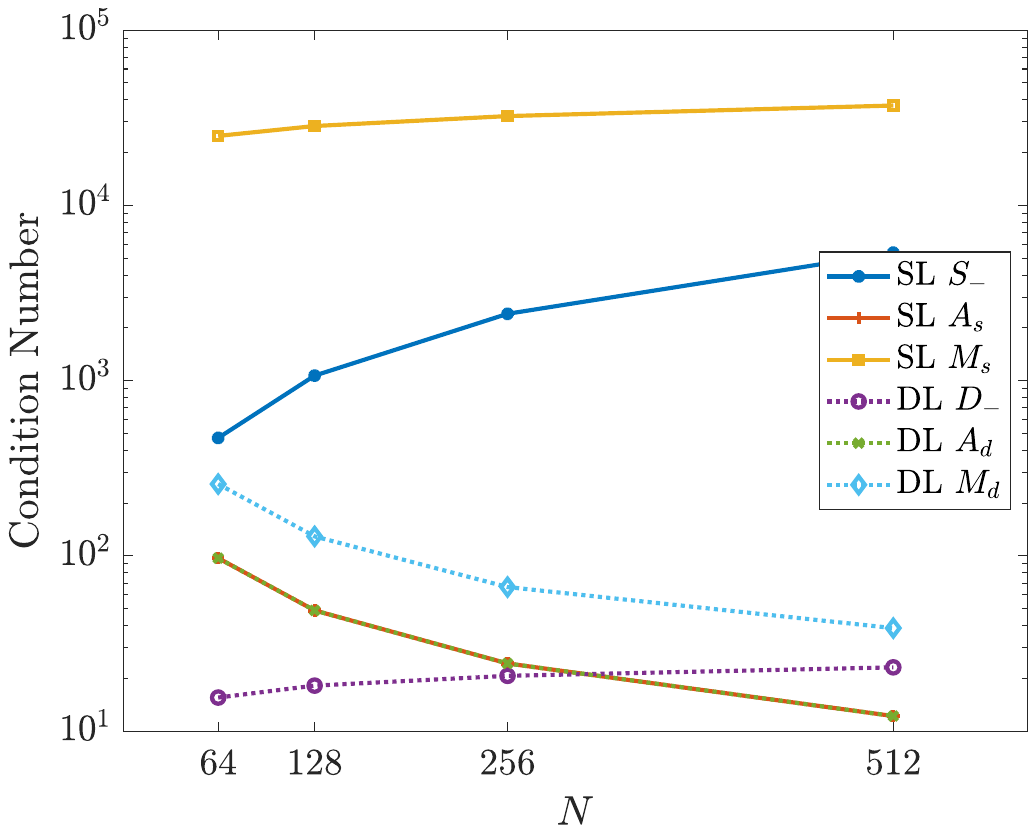}
        \caption{$\alpha=8$}
    \end{subfigure}
    \caption{Condition number (Dirichlet BC) for diamonds of different aspect ratios}\label{fig:diamond_cond}
\end{figure}

The condition numbers in Fig~\ref{fig:diamond_cond} are independent of the aspect ratios of the geometry, similar to the ellipse cases. The single layer formulations $S_-$ and $M_s$ have logarithmic growth, $A_s$ and double layer formulations $D_-, M_d$ and $A_d$ are uniformly bounded. In particular, $A_s, A_d$ and $M_d$ have logarithmic decay.

\begin{figure}[htbp]
\centering
    \begin{subfigure}{0.45\textwidth}
        \centering
        \includegraphics[width=\textwidth, trim = 2cm 7cm 2cm 6.5cm]{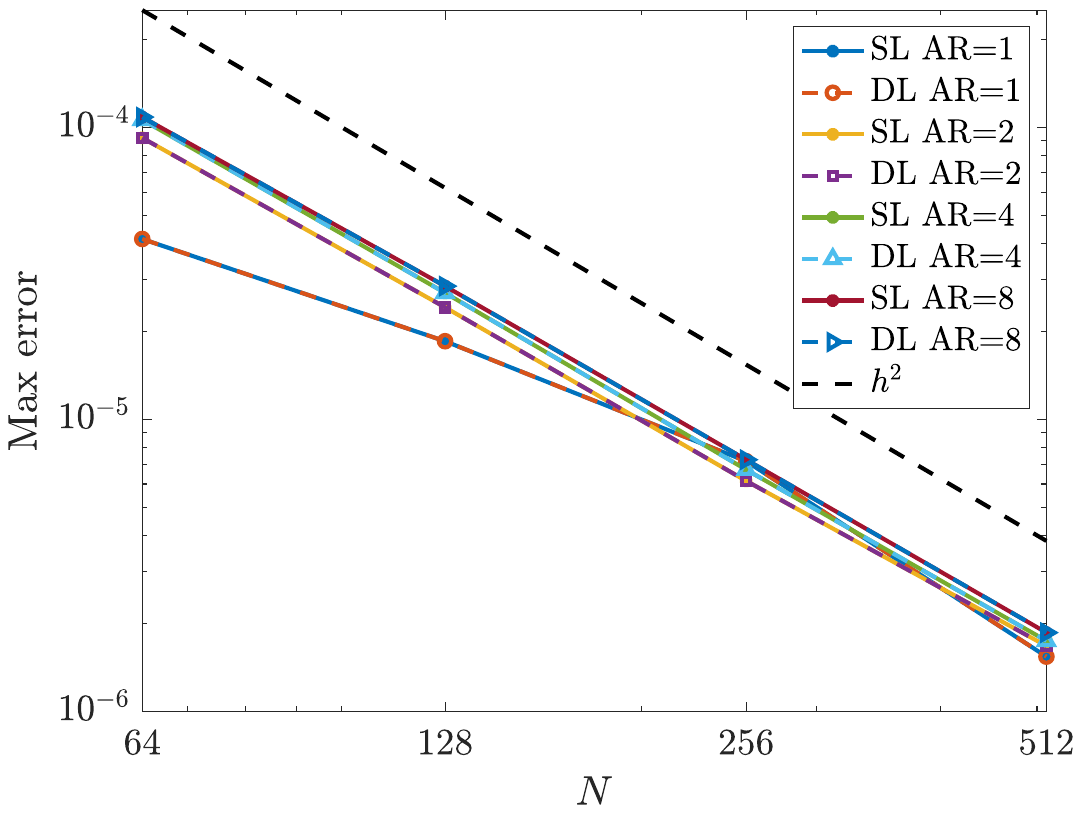}
        \caption{Convergence}\label{fig:diamond_conv}
    \end{subfigure}
    ~
    \begin{subfigure}{0.45\textwidth}
        \centering
        \includegraphics[width=\textwidth, trim = 2cm 7cm 2cm 6.5cm]{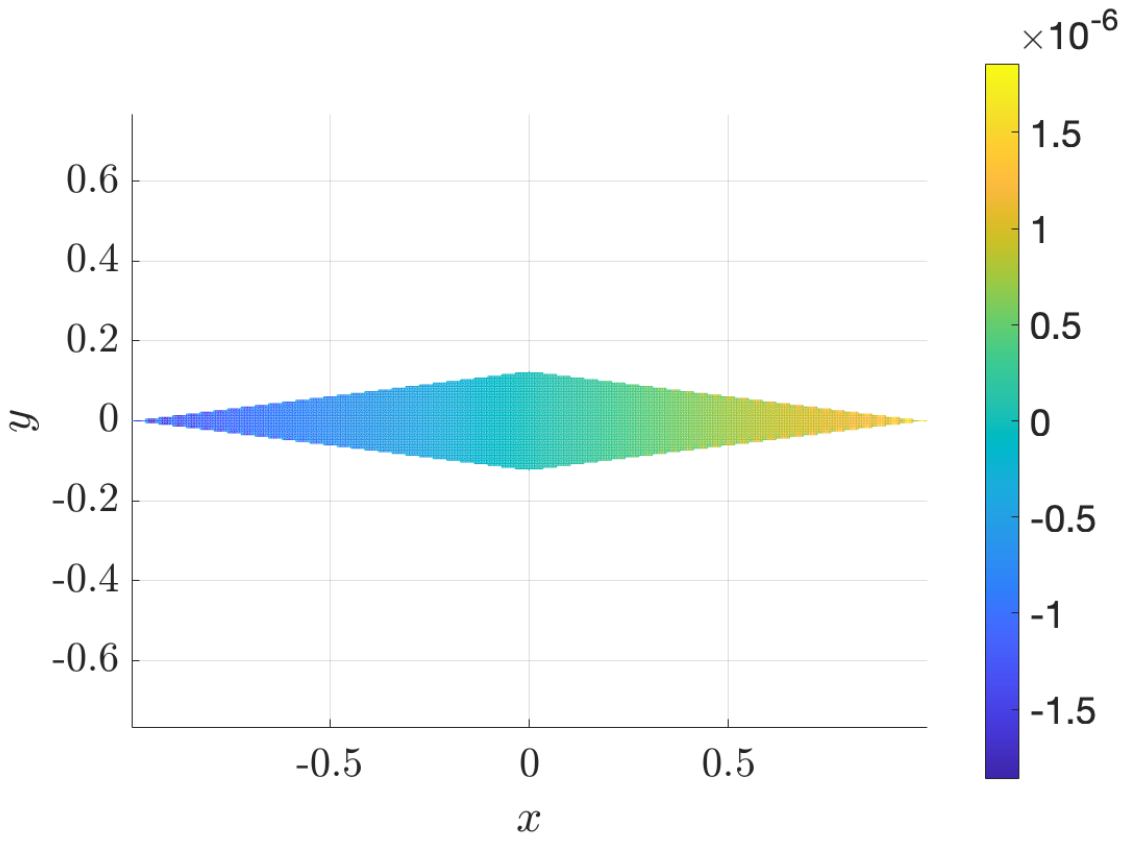}
        \caption{Errors for $\alpha=8$ ($512\times 512$)}\label{fig:diamond_error}
    \end{subfigure}
    \caption{Convergence (Dirichlet BC) for diamond of different aspect ratios and errors for aspect ratio $\alpha=8$}
\end{figure}

In Fig~\ref{fig:diamond_conv}, the second order convergence rate is also observed and the magnitudes of the max errors are also independent of the aspect ratios of the diamond or the single or double layer formulations, when the meshes are sufficiently refined. In Fig~\ref{fig:diamond_error}, the errors are also smooth, with larger absolute values occurring at where the solutions are large in magnitude. The sizes of the angles of the corners do not seem to affect the accuracy, as we increase the aspect ratios of the diamond.

\begin{remark}
Due to strong form the discretization of boundary conditions, the mesh does not see the corners. The developed method is yet to be tested when the solution itself exhibits singularities at corners such as in the Helmholtz equation.
\end{remark}





\subsection{Robin BC in ellipses}
Next, we test the Poisson equation with Robin BC
\begin{align}
\frac{\partial u(x)}{\partial n} + u(x) = g(x),\quad x\in \Gamma
\end{align} in ellipses of different aspect ratios.
\begin{figure}[htbp]
    \centering
    \begin{subfigure}{0.45\textwidth}
        \centering
        \includegraphics[width=\textwidth, trim = 2cm 7cm 2cm 6.5cm]{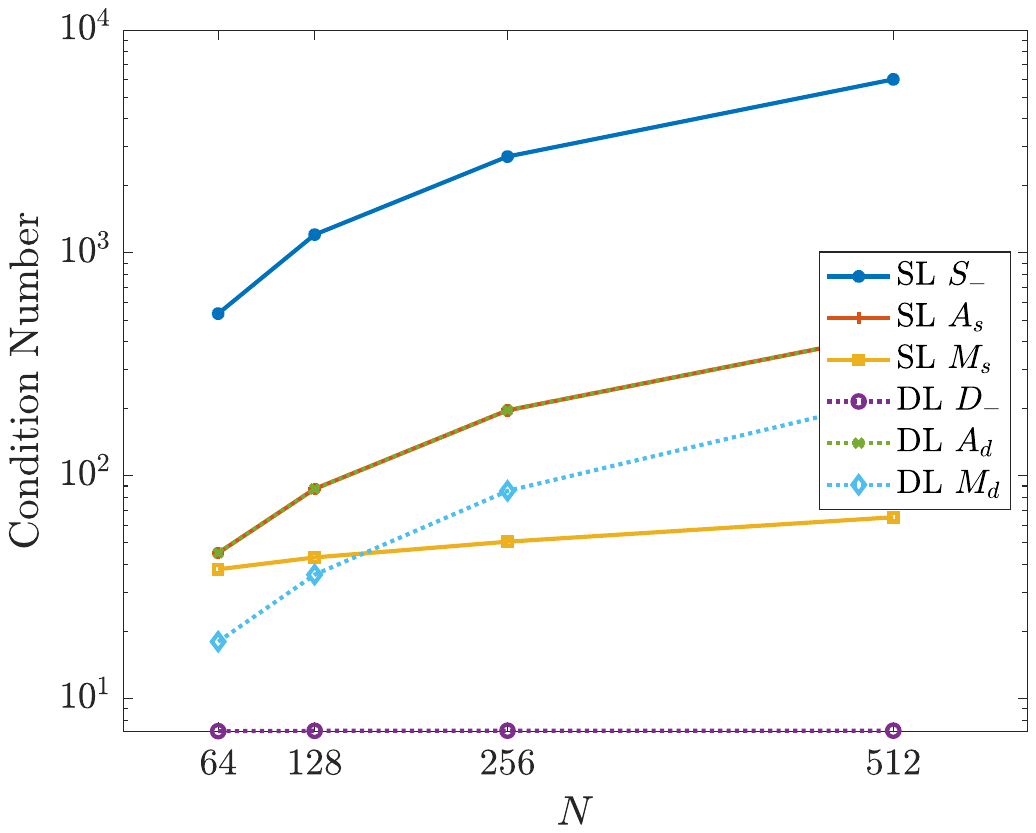}
        \caption{$\alpha=2$}
    \end{subfigure}
    ~
    \begin{subfigure}{0.45\textwidth}
        \centering
        \includegraphics[width=\textwidth, trim = 2cm 7cm 2cm 6.5cm]{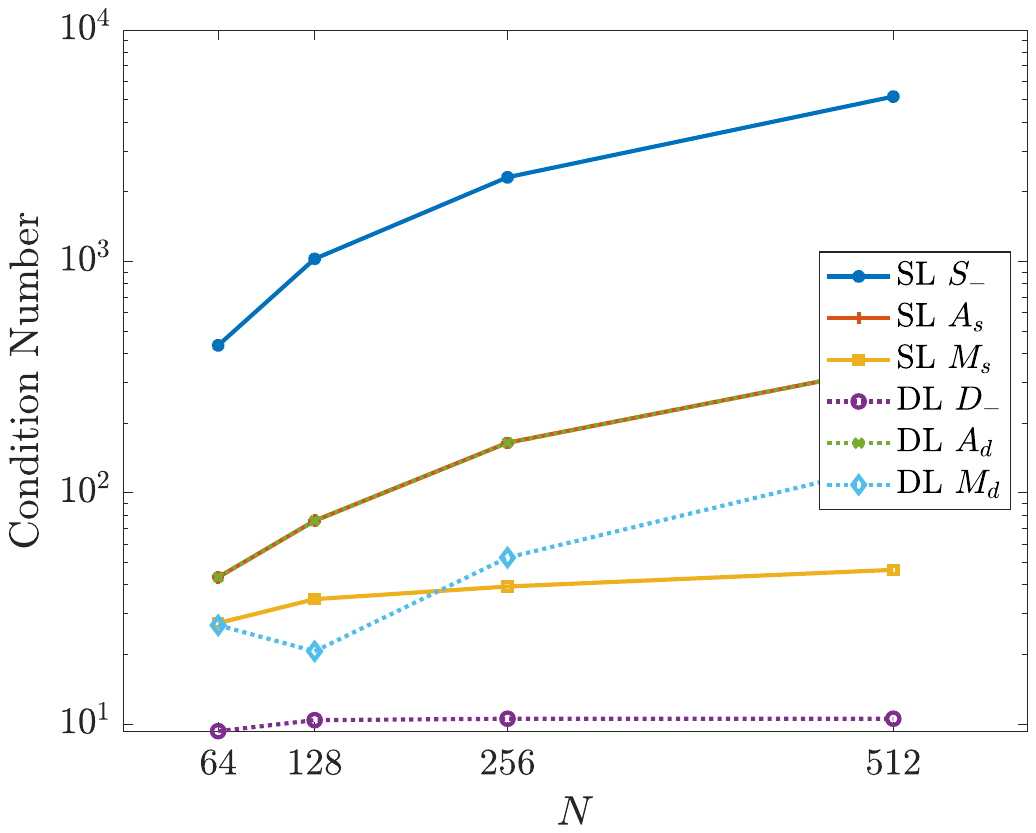}
        \caption{$\alpha=8$}
    \end{subfigure}
    \caption{Condition number for ellipses of different aspect ratios}\label{fig:robin_cond}
\end{figure}

Recall for Robin BC, the matrices are defined as 
\begin{align}
A_s&=\Phi_+S_+S^{-1}_{-} +\Phi_--\Phi'_-(R_+S_+S^{-1}_{-}+R_-),\\
M_s&=\Phi_+S_+ +\Phi_-S_{-}-\Phi'_-(R_+S_++R_-S_{-}),\\
A_d&=\Phi_+D_+D^{-1}_{-} +\Phi_--\Phi'_-(R_+D_+D^{-1}_{-}+R_-),\\
M_d&=\Phi_+D_+ +\Phi_-D_{-}-\Phi'_-(R_+D_++R_-D_{-}).
\end{align}
The patterns of the condition numbers in the Robin BC case (Fig~\ref{fig:robin_cond}) is different from the Dirichlet BC case. Only the double layer formulation $D_-$ has uniformly bounded condition number. The rest of the matrices $S_-,A_s,A_d,M_s,M_d$ follow logarithmic growth. The condition numbers are also independent of the aspect ratios of the ellipses.

\begin{figure}[htbp]
\centering
    \begin{subfigure}{0.45\textwidth}
        \centering
        \includegraphics[width=\textwidth, trim = 2cm 7cm 2cm 6.5cm]{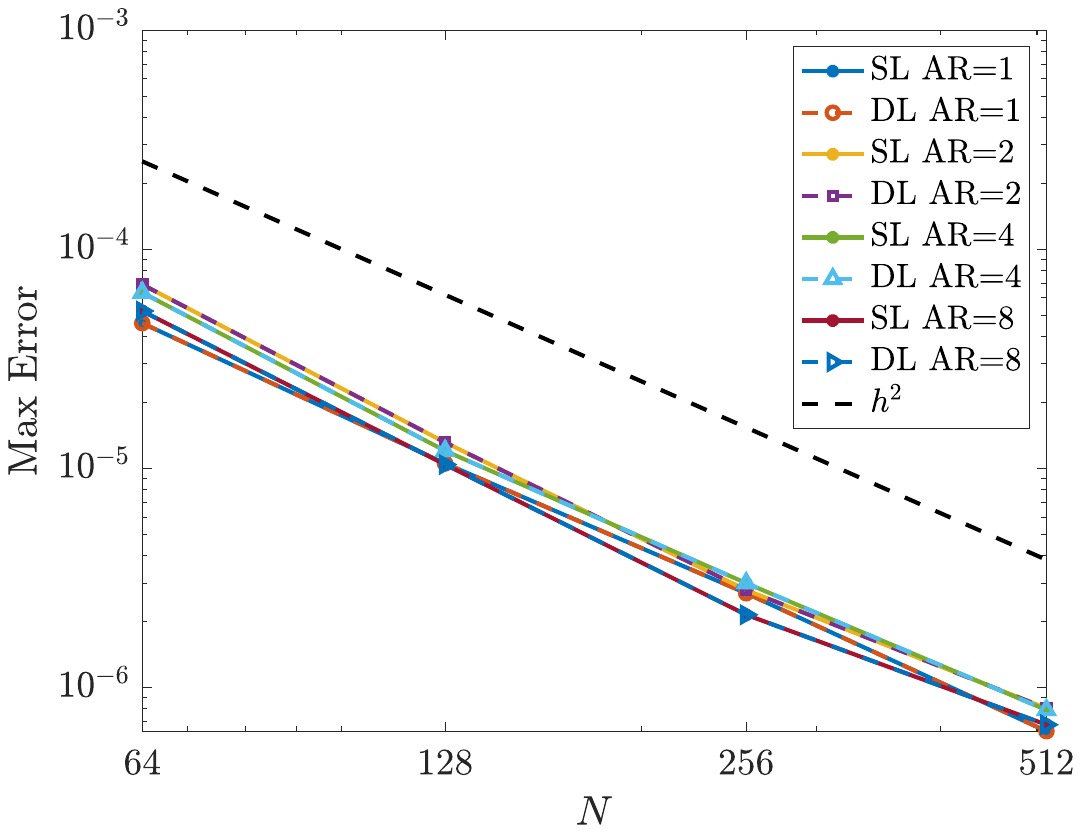}
        \caption{Convergence}\label{fig:robin_conv}
    \end{subfigure}
    ~
    \begin{subfigure}{0.45\textwidth}
        \centering
        \includegraphics[width=\textwidth, trim = 2cm 7cm 2cm 6.5cm]{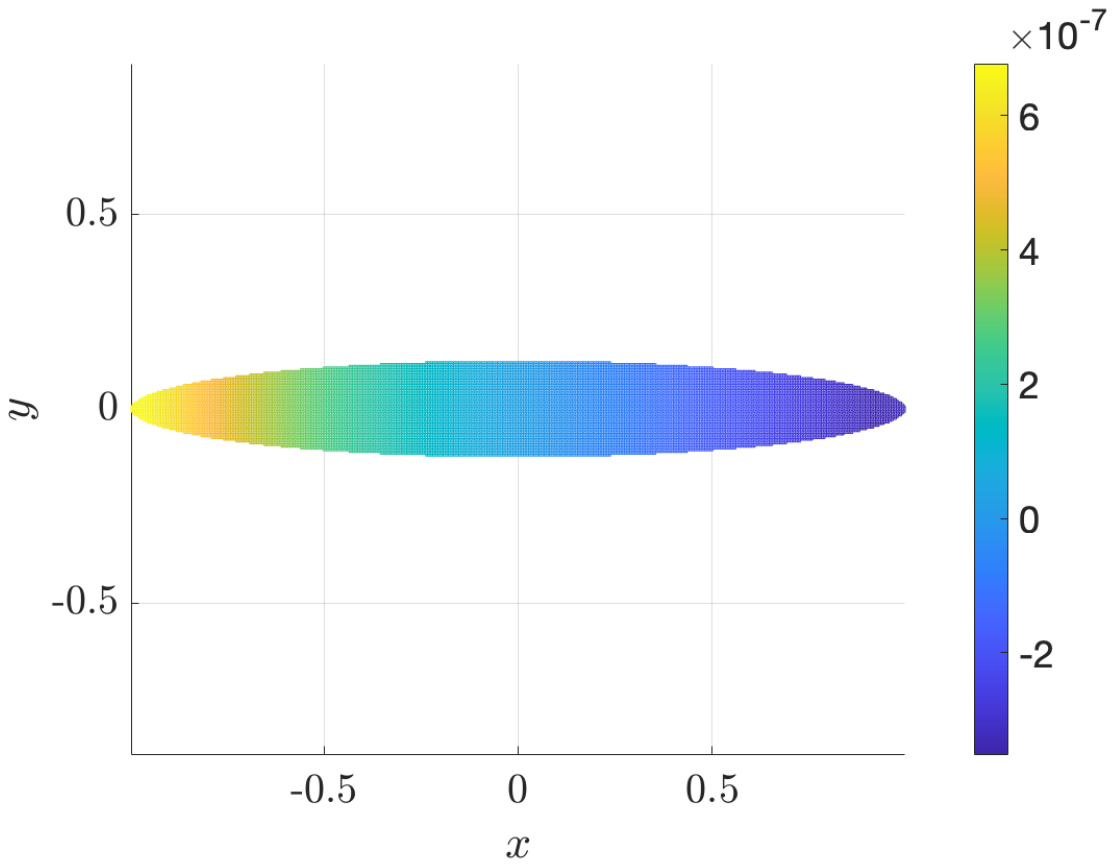}
        \caption{Errors for $\alpha=8$ ($512\times 512$)}\label{fig:robin_error}
    \end{subfigure}
    \caption{Convergence (Robin BC) for ellipse of different aspect ratios and errors for aspect ratio $\alpha=8$}\label{fig:robin}
\end{figure}
The convergence in Fig~\ref{fig:robin_conv} is standard second order. The magnitudes of the max errors also do not depend on the aspect ratios of the ellipses, regardless of whether single or double layer formulations are used. The errors in Fig~\ref{fig:robin_error} are also smooth with large errors occur near the tips where the solution itself is large in magnitude.


\subsection{Dirichlet and Neumann BC in unbounded domain}
Now we consider an irrotational and incompressible potential flow around a unit circle. One exact solution can be given as
\begin{align}\label{eqn:stream}
u(x,y) = \frac{x}{x^2+y^2},\quad x^2+y^2\geq1,
\end{align}
which satisfies the Laplace equation outside the unit circle and the far field boundary condition $u\rightarrow0$ as $r=\sqrt{x^2+y^2}\rightarrow\infty$.

The computational domain is selected to be $[-3,3]\times[-3,3]$.
No artificial boundary condition at $x,y=\pm3$ is needed, as the approximated solution is constructed directly from discrete densities defined at the exterior grid boundary $\gamma_-$. We test two cases with Dirichlet BC and Neumann BC respectively. The boundary conditions are computed using the exact solution \eqref{eqn:stream}, respectively.

We should also note that the double layer formulation fails in this case as the matrix $D_-$ will be singular.

\begin{figure}[htbp]
\centering
    \begin{subfigure}{0.45\textwidth}
        \centering
        \includegraphics[width=\textwidth, trim = 2cm 7cm 2cm 6.5cm]{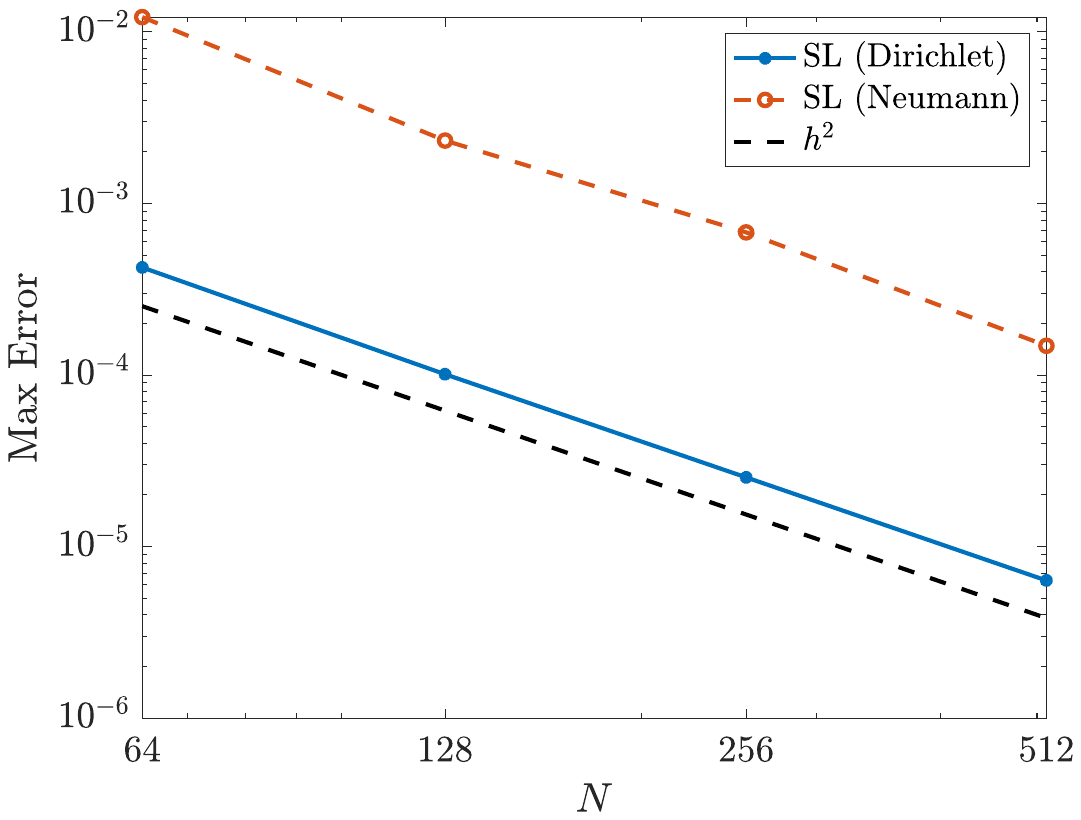}
        \caption{Convergence}\label{fig:unbounded_conv}
    \end{subfigure}
    ~
    \begin{subfigure}{0.45\textwidth}
        \centering
        \includegraphics[width=\textwidth, trim = 2cm 7cm 2cm 6.5cm]{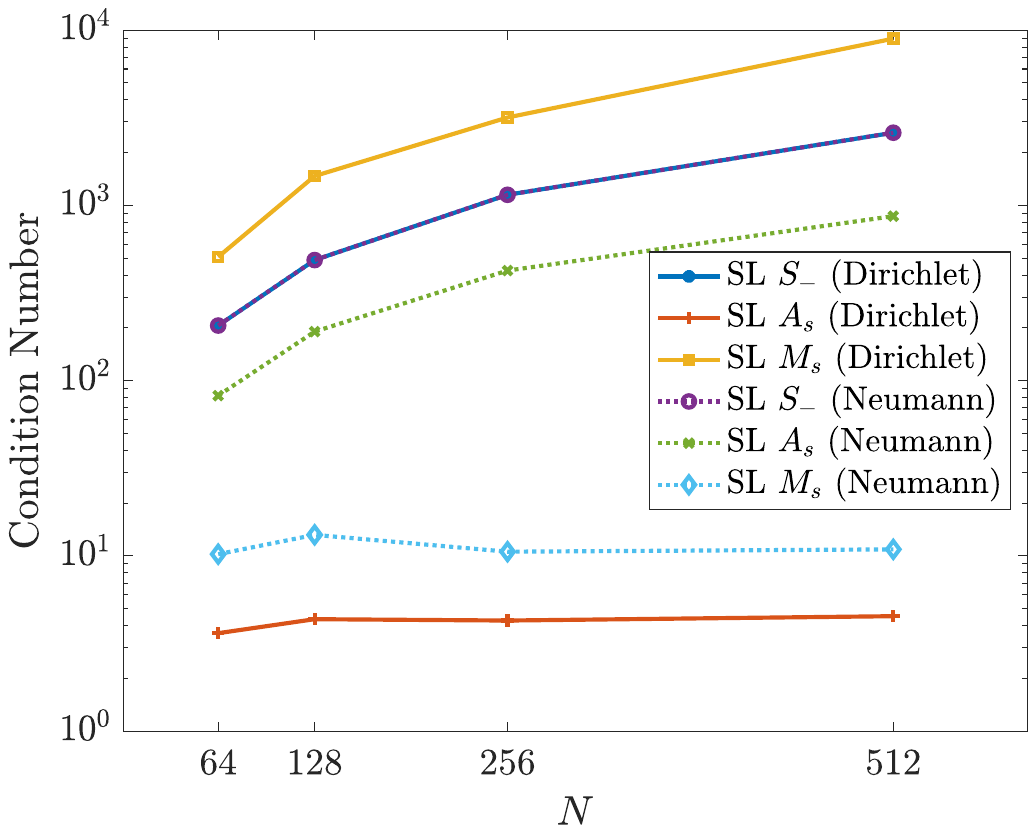}
        \caption{Condition number}\label{fig:unbounded_cond}
    \end{subfigure}
    \caption{Convergence and condition number for the unbounded domain}
\end{figure}

In Fig~\ref{fig:unbounded_conv}, we demonstrate the second order convergence rate of the max errors. The Neumann BC has slightly larger errors. In Fig~\ref{fig:unbounded_cond}, the matrices assembled for Dirichlet BC and Neumann BC exhibit different behaviors. All matrices have logarithmic growth except $A_s$ in the Dirichlet BC case (which is expected since we can regard it as right preconditioned) and $M_s$ in the Neumann BC case, which both are uniformly bounded.

\begin{figure}[htbp]
\centering
    \begin{subfigure}{0.45\textwidth}
        \centering
        \includegraphics[width=\textwidth, trim = 2cm 7cm 2cm 7.5cm]{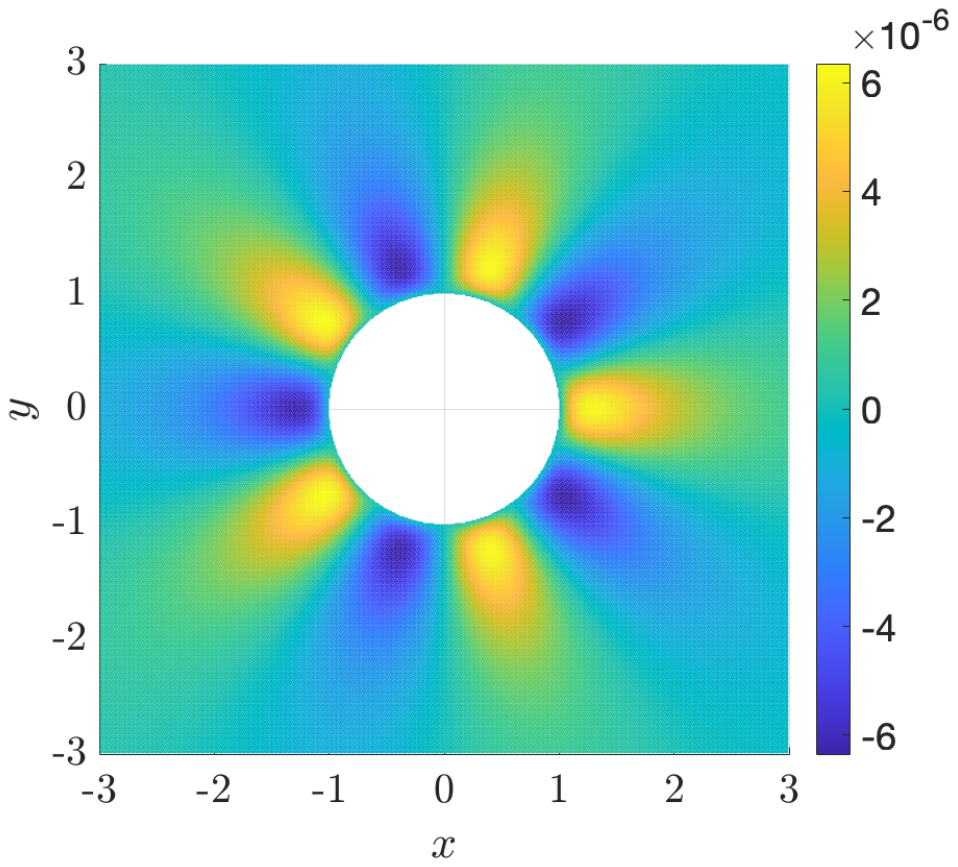}
        \caption{Dirichlet BC ($512\times 512$)}
    \end{subfigure}
    ~
    \begin{subfigure}{0.45\textwidth}
        \centering
        \includegraphics[width=\textwidth, trim = 2cm 7cm 2cm 7.5cm]{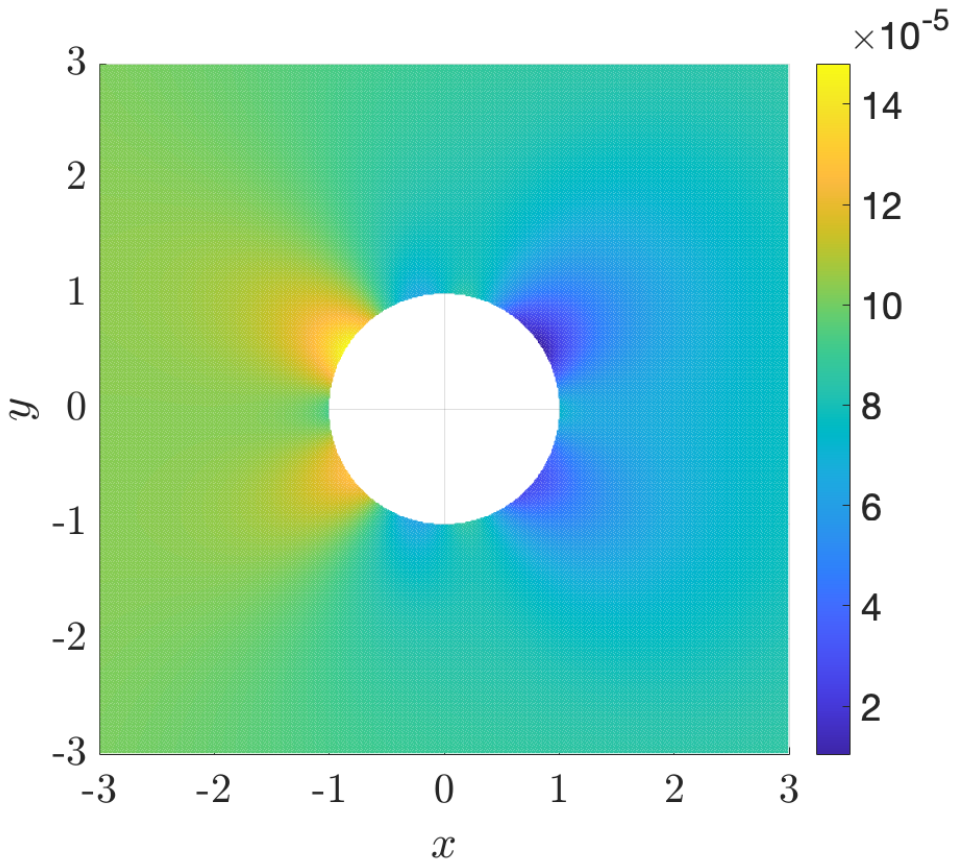}
        \caption{Neumann BC ($512\times 512$)}
    \end{subfigure}
    \caption{Errors for the unbounded domain}\label{fig:unbounded_errors}
\end{figure}
The error patterns are also different for Dirichlet BC and Neumann BC in Fig~\ref{fig:unbounded_errors}, where the errors from the Dirichlet BC exhibit periodicity along the circular boundary. Both cases seem to have the largest error close to the boundary.

\section{Conclusion}\label{sec:conclusion}
We have developed a computationally efficient unfitted boundary algebraic equation method that successfully extends the discrete potential theory framework to arbitrary complex geometries. The method combines lattice Green's functions with carefully constructed local basis functions on cut cells, achieving second-order accuracy for both Dirichlet, Neumann and Robin boundary conditions.

The theoretical analysis demonstrates that the method maintains stability and optimal convergence rates through a novel interpolation strategy that avoids extrapolation even near boundaries. In particular, the developed method does not suffer from small-cut cell issues as in the unfitted finite element method. The numerical result reveals mesh-independent conditioning for double layer formulations with Dirichlet conditions and controlled growth for Robin conditions, comparable to continuous boundary integral methods. In contrast, our unfitted method avoids evaluations of singular integrals

Numerical experiments confirm the theoretical predictions across diverse geometric configurations. The method handles ellipses and diamonds with extreme aspect ratios, sharp corners, and unbounded domains with uniform accuracy and without special numerical treatment. The error patterns remain smooth and predictable, with maximum errors independent of geometric aspect ratios or corner angles.

The integration with the difference potentials framework provides significant computational advantages. For bounded domains, the auxiliary problem formulation eliminates the need for direct evaluation of layer potentials at interior points, enabling the use of fast Fourier transform-based solvers. The method naturally extends to nonhomogeneous problems through method of superposition.

The unfitted BAE method fills a critical gap in numerical methods for elliptic problems, providing the geometric flexibility of unfitted methods with the favorable spectral properties of boundary integral formulations. Its robustness, efficiency, and ease of implementation make it particularly suitable for applications in computational physics and engineering where complex geometries are prevalent.

Future work includes extension of the current approach to more complicated applications, such as interface problems, biharmonic equations, Helmholtz equation and Stokes equation in bounded or unbounded domains. The matrix properties observed in the numerical results could also benefit from more rigorous spectral analysis.

\section*{Acknowledgement}

This work is partially funded by Natural Science Foundation of China (NSFC Grant No: 12401546) and Wenzhou Kean University (Grant No: ISRG2024003 and KY20250604000452).





\bibliographystyle{elsart-num}
\bibliography{references.bib}

\begin{thebibliography}{10}
\expandafter\ifx\csname url\endcsname\relax
  \def\url#1{\texttt{#1}}\fi
\expandafter\ifx\csname urlprefix\endcsname\relax\def\urlprefix{URL }\fi

\bibitem{peskin2002immersed}
C.~S. Peskin, The immersed boundary method, Acta numerica 11 (2002) 479--517.

\bibitem{leveque1994immersed}
R.~J. LeVeque, Z.~Li, The immersed interface method for elliptic equations with
  discontinuous coefficients and singular sources, SIAM Journal on Numerical
  Analysis 31~(4) (1994) 1019--1044.

\bibitem{hansbo2002unfitted}
A.~Hansbo, P.~Hansbo, An unfitted finite element method, based on {Nitsche}’s
  method, for elliptic interface problems, Computer methods in applied
  mechanics and engineering 191~(47-48) (2002) 5537--5552.

\bibitem{zhang2004immersed}
L.~Zhang, A.~Gerstenberger, X.~Wang, W.~K. Liu, Immersed finite element method,
  Computer Methods in Applied Mechanics and Engineering 193~(21-22) (2004)
  2051--2067.

\bibitem{ying2007kernel}
W.~Ying, C.~S. Henriquez, A kernel-free boundary integral method for elliptic
  boundary value problems, Journal of computational physics 227~(2) (2007)
  1046--1074.

\bibitem{glowinski1994fictitious}
R.~Glowinski, T.-W. Pan, J.~Periaux, A fictitious domain method for dirichlet
  problem and applications, Computer Methods in Applied Mechanics and
  Engineering 111~(3-4) (1994) 283--303.

\bibitem{helsing2013solving}
J.~Helsing, Solving integral equations on piecewise smooth boundaries using the
  rcip method: a tutorial, in: Abstract and applied analysis, Vol. 2013, Wiley
  Online Library, 2013, p. 938167.

\bibitem{klockner2013quadrature}
A.~Kl{\"o}ckner, A.~Barnett, L.~Greengard, M.~O'Neil, Quadrature by expansion:
  A new method for the evaluation of layer potentials, Journal of Computational
  Physics 252 (2013) 332--349.

\bibitem{duffin1953discrete}
R.~J. Duffin, Discrete potential theory, Duke Math. J. 20~(1) (1953) 233--251.

\bibitem{saltzer1958discrete}
C.~Saltzer, Discrete potential theory for two-dimensional {Laplace} and
  {Poisson} difference equations, Tech. rep. (1958).

\bibitem{tsynkov2003definition}
S.~V. Tsynkov, On the definition of surface potentials for finite-difference
  operators, Journal of scientific computing 18 (2003) 155--189.

\bibitem{martinsson2009boundary}
P.-G. Martinsson, G.~J. Rodin, Boundary algebraic equations for lattice
  problems, Proceedings of the Royal Society A: Mathematical, Physical and
  Engineering Sciences 465~(2108) (2009) 2489--2503.

\bibitem{economou2006green}
E.~N. Economou, Green's functions in quantum physics, Vol.~7, Springer Science
  \& Business Media, 2006.

\bibitem{pozrikidis2014introduction}
C.~Pozrikidis, An introduction to grids, graphs, and networks, Academic, 2014.

\bibitem{ryaben2012method}
V.~S. Ryaben'Kii, Method of {Difference} {Potentials} and its applications,
  Vol.~30, Springer Science \& Business Media, 2012.

\bibitem{larsson2003partial}
S.~Larsson, V.~Thom{\'e}e, Partial differential equations with numerical
  methods, Vol.~45, Springer, 2003.

\bibitem{martinsson2002asymptotic}
P.-G. Martinsson, G.~J. Rodin, Asymptotic expansions of lattice green's
  functions, Proceedings of the Royal Society of London. Series A:
  Mathematical, Physical and Engineering Sciences 458~(2027) (2002) 2609--2622.

\bibitem{cserti2000application}
J.~Cserti, Application of the lattice {Green}’s function for calculating the
  resistance of an infinite network of resistors, American Journal of Physics
  68~(10) (2000) 896--906.

\bibitem{morita1971useful}
T.~Morita, Useful procedure for computing the lattice green's function-square,
  tetragonal, and bcc lattices, Journal of mathematical physics 12~(8) (1971)
  1744--1747.

\bibitem{gillman2014fast}
A.~Gillman, P.-G. Martinsson, A fast solver for {Poisson problems} on infinite
  regular lattices, Journal of Computational and Applied Mathematics 258 (2014)
  42--56.

\bibitem{gabbard2024lattice}
J.~Gabbard, W.~M. van Rees, Lattice {Green}’s functions for high-order finite
  difference stencils, SIAM Journal on Numerical Analysis 62~(1) (2024) 25--47.

\bibitem{xia2023local}
Q.~Xia, Local-basis difference potentials method for elliptic {PDEs} in complex
  geometry, Journal of Computational Physics 488 (2023) 112246.

\bibitem{banks2016galerkin}
J.~W. Banks, T.~Hagstrom, On {Galerkin} difference methods, Journal of
  Computational Physics 313 (2016) 310--327.

\bibitem{helsingfast2011}
J.~Helsing, A fast and stable solver for singular integral equations on
  piecewise smooth curves, SIAM Journal on Scientific Computing 33~(1) (2011)
  153--174.

\bibitem{bremer2010efficient}
J.~Bremer, V.~Rokhlin, Efficient discretization of laplace boundary integral
  equations on polygonal domains, Journal of Computational Physics 229~(7)
  (2010) 2507--2525.

\end{thebibliography}
 
\end{document}